\documentclass[12pt,a4paper]{article}
\usepackage{lineno,hyperref}
\usepackage[latin5]{inputenc}
\usepackage{amsmath}
\usepackage{amssymb}
\usepackage{graphicx}
\usepackage{subfigure}
\usepackage{caption}
\usepackage{float}
\usepackage{mathtools}
\usepackage{color}
\usepackage{cite}
\usepackage{mathrsfs}
\usepackage[]{geometry}
\usepackage{subfigure}
\DeclareMathOperator{\sech}{sech}
\title{Stability Satisfied Numerical Approximates to the Non-analytical Solutions of the Cubic Schrödinger Equation}
\author{Alper Korkmaz\thanks{alperkorkmaz7@gmail.com} \\%$^{,}$\thanks{Senior Author}\\
{\scriptsize Çankırı Karatekin University, Department of Mathematics, 18200, Çankırı, Turkey.}} 
\begin{document}
\maketitle
\begin{abstract}

\noindent
The time dependent complex Schrödinger equation with cubic nonlinearity is solved by constructing differential quadrature algorithm based on sinc functions. Reduction to a coupled system of real equations enables to approach the space derivative terms by the proposed method. The resulted ordinary differential equation system is integrated with respect to the time variable by using  a bunch explicit methods of lower and higher orders. Some initial boundary value problems containing some analytical and non-analytical initial data are solved for experimental illustrations. The computational errors between the analytical and numerical solutions are measured by the discrete maximum error norm in case the analytical solution exists. The two conserved quantities are calculated by using the numerical results in all cases. The matrix stability analysis is implemented to control the time step size. 

\end{abstract}

\section{Introduction}
\noindent
Consider the cubic nonlinear Schrödinger (NLS) equation of the canonical form \cite{zakharov1}
\begin{equation}
i\frac{\partial u(x,t)}{\partial t}+\frac{\partial ^2 u(x,t)}{\partial x^2}+\kappa \left| u(x,t)\right|^2u(x,t)=0, \, -\infty<x<\infty, \, t>0\label{nls}
\end{equation}
where $u=u(x,t) \in \mathbb{C}$. This equation is a standard dimensionless form of the Talanov's model \cite{talanov} describing stationary light beam of plane form in a medium with nonlinear refractive index and quasimonochromatic one-dimensional wave in a dispersive and inertialess nonlinear medium\cite{zakharov1,bespalov,zakharov2}. The solutions are named as self-focusing or defocussing according to the sign of the coefficient of the nonlinear term $\kappa$ in the propagating of electromagnetic waves\cite{peregrine1}. The NLS equation disallows steady soliton solutions traveling with a constant speed\cite{debnath1}. Alternatively, it has envelop soliton solutions containing both exponential and $\sech$ functions traveling with different speeds\cite{debnath1}. Zakharov and Shabat\cite{zakharov1} applied the inverse scattering method successfully to develop exact solutions of some dynamic problems to the cubic NLS equation (\ref{nls}). Debnath\cite{debnath1} summarizes the findings of Zakharov and Shabat\cite{zakharov1} as
\begin{itemize}
\item An initial envelope pulse breaks into permanent, propagating, shorter scaled solitons the number of which depends on the initial data and an oscillatory tail,
\item Separates after interaction without changes except possible position or phase changes,
\item The oscillatory tail whose nature is defined in the initial data disperses linearly with a decaying amplitude as $t \rightarrow \infty$,
\end{itemize}
when the initial data converge zero as $|x|\rightarrow\infty$.

\noindent
An envelop-type soliton also behaves like a particle familiar to the KdV soliton\cite{yajima1}. The amplitude of the envelop of the one soliton solution of the cubic NLS equation (\ref{nls}) is of permanent form describing a physical wave which is in a good agreement with Hammack's unpublished experiment data\cite{ablowitz1}.

\noindent
Recurrence is also significant for the solutions of the NLS equation in bounded or periodic domains\cite{debnath1}. The relation between the recurrence, dimensionality, and the stability in the Lagrange sense of solutions of the cubic NLS equation is investigated in details by Thyagaraja\cite{thy1}.

\noindent
The equation can be a model for the vortex line in an imperfect Bose gas with weak pair repulsions between atoms\cite{pita1}. In addition to be an approximation for beams of the modulated form in nonlinear optics, the NLS is a significant equation for the time dependent dispersive waves\cite{whitham1}.

\noindent
Tsuziki\cite{tsuziki1} examines deeply the nonlinear solitary and periodic waves for a particular type of the NLS equation. Some conservations laws describing the quantities density, current and energy are derived in integral form for that particular type. The collision of two positive solitaries with different heights traveling to the right in the horizontal axis is also discussed. In the last part of this study, the wave generation in terms of the decay of an arbitrary initial disturbance into solitary waves is investigated, too.

\noindent
Some of the exact solutions covering interaction of solitary wave packets and $n$-soliton solutions of the cubic NLS equation are derived by transforming it to a problem of inverse scattering. Some conservation laws describing various quantities for the NLS equation (\ref{nls}) are also defined in \cite{zakharov1}. The direct scattering problem for the NLS equation (\ref{nls}) is solved for a particular initial data $u(x,0)=\sech{x}$ with the condition $\kappa \geq 2$\cite{miles1}. A Cauchy problem constructed on the generalized function type initial data is solved in a certain algebra of some particular generalized functions\cite{bu1}. 

\noindent
Solitons of self-focusing form of the NLS equation are in different forms with various characteristic properties. The envelope soliton decays zero as a Ma soliton decays to the uniform solution\cite{peregrine1}. When the real components of the eigenvalues are equal, bi-soliton type solution whose velocity depends on these eigenvalues can be derived. On the other hand, a dark soliton decays to a uniform solution of the defocussing NLS equation\cite{zakharov3}. Peregrine's study\cite{peregrine1} also suggests a new type analytical solution in rational form describing an isolated amplitude peak. Some exact solutions expressed in terms of some trigonometric and hyperbolic functions can be constructed by He's Lindstedt-Poincar\'{e} method in the modified form\cite{ozis1}. The complex tangent function is capable to find some exact traveling wave solutions in the hyperbolic function form as some ansatzes defined as the multiplication of $\tanh$ and exponential functions leads some explicit exact solutions\cite{khuri1}.

\noindent
In addition to many analytical or theoretical studies on the NLS equation, various numerical methods have been derived for the analytical or non-analytical solutions.

\noindent
Taha proposes some local and global methods for the numerical solution to an initial boundary value problem whose solution is a model for the motion of an single initial pulse. He compares his results with the results obtained by some classical numerical methods covering local scheme, pseudospectral and split step Fourier methods\cite{taha1}.

\noindent
Dereli et al.\cite{dereli1} investigate the motion of single positive soliton and the collision of a couple of positive solitons moving in the opposite directions along the horizontal axis numerically. In the study, the radial basis meshless collocation method in four different radial functions, Gaussian, multiquadric, inverse quadric and inverse multiquadric is implemented. The authors report that the least error, accordingly, the best results are generated by the Gaussian functions.

\noindent
Some finite difference methods covering linearized Crank-Nicolson scheme are derived for solutions of the three model problems for the inhomogeneous NLS equation by Chang et al.\cite{chang1}. A comparison with Hopscotch-type methods, split step Fourier and spectral schemes indicates that the newly proposed method is efficient and robust.   

\noindent
Quadratic B-spline finite element method\cite{dag1} is another significant study in the related literature. Dag obtains the numerical simulations of various initial boundary value problems covering propagation of single soliton, bound state of solitons and wave birth models for both standing and traveling initial pulses. This paper is one of the earlier numerical studies dealing with the NLS equation. Gardner et al.\cite{gardner1} solve different initial-boundary value analytical and non-analytical problems with various characteristics for the cubic NLS equation by the cubic B-spline finite elements. The Taylor collocation method based on the quintic B-splines is derived for the similar problems of the cubic NLS equation\cite{aksoy1}. 

\noindent
Bound states of the NLS equation are deeply investigated by using $L^2$-Galerkin with product approximation and Ablowitz \& Ladik integrable finite difference methods in \cite{robinson1}. Twizell et al.\cite{twi1} reduce some initial-boundary value problems for the nonlinear cubic NLS equation to a linear initial value problems of order one by a family of finite difference techniques. They also examined the truncation error, stability and convergence properties of the proposed methods in details. Moreover, they suggest a report explaining the effects of the coefficients to the error between the numerical and the analytical solutions.

\noindent
Some wave birth models with small dispersion parameters to the focusing NLS equation is studied in details with a second order semi-implicit adaptive moving mesh method\cite{hector1}.

\noindent
Differential quadrature methods based on cosine expansion have been successfully implemented to five analytical or non-analytical problems for the NLS equation\cite{korkmaz1}. Motion of single soliton, collision of two positive solitons traveling in the opposite directions along the horizontal axis, bound state of solitons, wave birth by standing or traveling single solitary wave are simulated successfully. The conservation laws are in a good agreement with the theoretical aspects as expected. In the same year, Korkmaz and Dag\cite{korkmaz2} announce that they solve the same equation with a variation of the differential quadrature method. This time, they use Lagrange interpolation polynomials as basis in the space discretization but do not change the time integration technique. %These two studies also includes comparison of implemented methods with some earlier numerical studies.

\noindent
Different from the last two studies, we derive differential quadrature method based on sine cardinal functions combined with a family of time integration techniques in different classes for the numerical solutions of some initial boundary value problems to cubic NLS equation. The assumption $u(x,t)=f(x,t)+ig(x,t),\, i=\sqrt{-1}$ reduces the cubic NLS equation (\ref{nls}) to a coupled system of ODEs 
\begin{equation}
\begin{aligned}
g_t&=f_{xx}+\kappa (f^2+g^2)f \\
f_t&=-g_{xx}-\kappa (f^2+g^2)g
\end{aligned}\label{nlss}
\end{equation}
where $f=f(x,t)$ and $g=g(x,t)$ are real functions. The artificial homogeneous Dirichlet boundary data at both ends are completely compatible to the chosen models having physical requirement that $u(x,t)\rightarrow 0$ as $|x|\rightarrow \infty$. Since $u(x,t)=f(x,t)+ig(x,t)$, the boundary conditions for the system (\ref{nlss}) are adapted as
\begin{equation}
\begin{aligned} 
f(a,t)&=0,\, f(b,t)=0 \\
g(a,t)&=0,\, g(b,t)=0 \\
\end{aligned} \label{bc}
\end{equation}

\section{Design of the Method}

\noindent
Consider a sufficiently smooth function $u=u(x,t)$ over a finite interval $[a,b]$. Even though the function $u$ is function of two variables, the second variable $t$ (the time variable through the study) is assumed to be fixed while approximating to the derivatives with respect to the space variable $x$. Thus, $r.$th order the derivative of the function $u$ with respect to $x$ is approximated by the finite weighted sum of all functional values in the interval $[a,b]$. Let $[a,b]$ be partitioned as $P:a=x_0<x_1<\ldots<x_N=b$. One should mention that each grid can be written in terms of grid size and the subscript of the grid as $x_m=m\Delta x$. %The partition $P$ is not necessary to be uniform. 
The definition of the differential quadrature derivative approximation of $u(x,t)$ at a distinct grid point $x_m$ is
\begin{equation}
\left. \dfrac{\partial ^{r}u(x,t)}{\partial x^{r}}\right|_{x=x_m}=\sum\limits_{j=0}^{N}{w^{(r)}_{m,j}u(x_j,t)},\, 0\leq m\leq N \label{dqm}
\end{equation} 
where $w^{(r)}_{m,j}$ is the weight of $u(x_j,t)$ for the $r.$th order derivative approximation of the function $u(x,t)$ at the internal grid point $x_i$\cite{belman1}. The significant point of the approximation is the determination of the weights $w^{(r)}_{m,j}$. Once, they are determined, the approximation (\ref{dqm}) is directly substituted instead of the related derivatives in the differential equation. Since there may exist different basis function sets spanning the same function or vector space, those basis all enable to determine the weights\cite{belman1,dqm1,dqm2,dqm3,dqm4,dqm5}. In this study, the weights are calculated by the set of sine cardinal functions spanning the problem interval.

\noindent
A sine cardinal function set $\{ T_m(x)\}_{m=0}^{N}$ with elements 
\begin{equation}
T_{m}(x)=\left\{ 
							\begin{array}{lcc}
										\dfrac{\sin{([\dfrac{x-m\Delta x}{\Delta x}]\pi)}}{[\dfrac{x-m\Delta x}{\Delta x}]\pi} & , & x \neq m\Delta x \\ 
										1 & , & x=m\Delta x \\ 
							\end{array}%
					\right.  \label{sinc}
\end{equation}
constitutes a basis for the functions defined in $[a,b]$ where $\Delta x$ is the equal grid size \cite{stenger,carlson1,secer1,dehghan2}. A sine cardinal function value at a grid in $[a,b]$ can be calculated easily as
\begin{equation}
T_m(x_j)=\delta_{mj} %= \left\{ 							\begin{array}{lcc}										1& , & m=j\\ 										0 & , & m \neq j \\ 							\end{array}%					\right.  
\label{sincnodal}
\end{equation}
where $\delta_{mj}$ is the Kronecker operator\cite{dehghan2}.
\noindent
The function $C(u)(x)$ approximated by an infinite convergent series
\begin{equation}
C(u)(x)= \sum_{m=-\infty}^{\infty}u(m\Delta x)T_m(x)        
\label{cardinal}
\end{equation} 
is named the cardinal of $u$ on ($-\infty ,\infty$) and it interpolates $u$ at the points that are integer multiple of $\Delta x$\cite{lund1}.

\noindent
The lowest ordered two derivatives of a sine cardinal function $T_m(x)$ are calculated in an explicit form as:
\begin{equation}
T_{m}^{\prime}(x)=\left\{ 
							\begin{array}{lcc}
										\dfrac{\dfrac{\pi}{\Delta x}(x-m\Delta x)\cos{\dfrac{x-m\Delta x}{\Delta x} \pi}-\sin{\dfrac{x-m\Delta x}{\Delta x} \pi} }{\dfrac{\pi}{\Delta x}(x-m\Delta x)^2}& , & x \neq m\Delta x \\ 
										0 & , & x=m\Delta x \\ 
							\end{array}%
					\right.  \label{sincd}
\end{equation}
\begin{equation}
T_{m}^{\prime \prime}(x)=\left\{ 
							\begin{array}{lcc}
							     \dfrac{-\dfrac{\pi}{\Delta x}\sin{\dfrac{x-m\Delta x}{\Delta x} \pi}}{x-m\Delta x}-\dfrac{2\cos{\dfrac{x-m\Delta x}{\Delta x} \pi}}{(x-m\Delta x)^2}+\dfrac{2\sin{\dfrac{x-m\Delta x}{\Delta x} \pi}}{\dfrac{\pi}{\Delta x}(x-m\Delta x)^3}& , & x \neq m\Delta x \\ 
									-\dfrac{{\pi}^2}{3\Delta x^2} & , & x=m\Delta x \\ 
							\end{array}%
					\right.  \label{sincdd}
\end{equation}

\noindent
In order to calculate the weights $w_{mj}^{(2)}$ of the approximation to second order derivative terms $f_{xx}$ and $g_{xx}$, we substitute each sine cardinal function into the differential quadrature approximation equation (\ref{dqm}) for $r=2$. 

\noindent
Assume that $m=0$. That means all $N+1$ weights $w_{0j}^{2}$ are determined by substituting each sine cardinal function $T_m(x)$ and its second order derivative into (\ref{dqm}). Although its not necessary to follow an order, all calculations are completed in an order for convenience and simplicity. Substituting $T_0(x)$ and its second order derivative into the (\ref{dqm}) gives
\begin{equation}
\begin{aligned}
T_0^{\prime \prime}(x_0)&=\sum \limits_{j=0}^{N}w_{0j}^{(2)}T_0(x_{j}) \\
     						 &=w_{00}^{(2)}{T_0(x_{0})}+w_{01}^{(2)}{T_0(x_{1})}+\ldots +w_{0N}^{(2)}{T_0(x_{N})} \\
								 &=w_{00}^{(2)}\delta_{00}+w_{01}^{(2)}\delta_{01}+\ldots +w_{0N}^{(2)}\delta_{0N} \\
							   \dfrac{-{\pi}^2}{3\Delta x^2}&=w_{00}^{(2)} 		 
\end{aligned}
\label{b01}
\end{equation}

\noindent
The same methodology can be used to find the weight $w_{01}^{(2)}$ as
\begin{equation}
\begin{aligned}
T_1^{\prime \prime}(x_0)&=\sum \limits_{j=0}^{N}w_{0j}^{(2)}T_0(x_{j}) \\
     						              &=w_{00}^{(2)}{T_1(x_{0})}+w_{02}^{(2)}{T_1(x_{1})}+\ldots +w_{0N}^{(2)}{T_1(x_{N})} \\
								              &=w_{00}^{(2)}\delta_{10}+w_{01}^{(2)}\delta_{11}+\ldots +w_{0N}^{(2)}\delta_{1N} \\
-\dfrac{2\cos{\left(\dfrac{x_0-\Delta x}{\Delta x}\pi\right)}}{(x_0-\Delta x)^2}&=w_{01}^{(2)} \\
-\dfrac{2\cos{\left(\dfrac{0\Delta x-\Delta x}{\Delta x}\pi\right)}}{(0\Delta x-\Delta x)^2}&=w_{01}^{(2)} \\ 
-\dfrac{2\cos{\left((0-1)\pi\right)}}{(0-1)^2 \Delta x^2}&=w_{01}^{(2)}\\
\dfrac{2(-1)^{(0-1+1)}}{(0-1)^2 \Delta x^2}&=w_{01}^{(2)} 		 	 
\end{aligned}
\label{b02}
\end{equation}
where the identity $\cos{(k\pi)}=(-1)^k$ is used to rearrange the left hand side of the equation. This procedure can generalized for any weight $w_{0m}^{(2)}$ related to the point $x_0$ originated from the point $x_m$ by using the basis $T_m(x)$ as
\begin{equation}
\begin{aligned}
T_m^{\prime \prime}(x_0)&=\sum \limits_{j=0}^{N}w_{0j}^{(2)}T_m(x_{j}) \\
     						              &=w_{00}^{(2)}{T_m(x_{0})}+w_{01}^{(2)}{T_m(x_{1})}+\ldots+w_{0m}^{(2)}{T_m(x_{m})}+\ldots +w_{0N}^{(2)}{T_m(x_{N})} \\
								              &=w_{00}^{(2)}\delta_{m0}+w_{01}^{(2)}\delta_{m1}+\ldots +w_{0m}^{(2)}\delta_{mm}+\ldots +w_{0N}^{(2)}\delta_{mN} \\
\dfrac{2(-1)^{(0-m+1)}}{(0-m)^2 \Delta x^2}&=w_{0m}^{(2)} 		 	 
\end{aligned}
\label{bm2}
\end{equation}
This explicit form can be extended for an arbitrary weight $w_{mj}^{(2)}$, related to the point $x_m$ originated from $x_j$, as
\begin{equation}
w_{mj}^{(2)}=\dfrac{2(-1)^{m-j+1}}{\Delta x^{2}(m-j)^{2}}
\end{equation} 
when $m\neq j$, and 
\begin{equation}
w_{mm}^{(2)}=-\dfrac{\pi ^{2}}{3\Delta x^{2}}
\end{equation}
when $m=j$. This explicit formulation of the weights are used some earlier studies for various problems \cite{bellomo1,korkmaz4}.

\noindent
Even though the determination of only the second derivatives weights is sufficient due to the structure of the NLS equation, the weights $w_{mj}^{(1)}$ should also be calculated to compute the conservation law $C_3$ having first order derivative in. Thus, we start by letting $r=1$ in the (\ref{dqm}). In order to determine the weights $w_{0j}^{(1)}$, assume that $m=0$ initially. Substituting the first element $T_0(x)$ of the basis functions set into the differential quadrature approximation (\ref{dqm}) gives
\begin{equation*}
\begin{aligned}
T_0^{\prime}(x_0)&=\sum \limits_{j=0}^{N}w_{0j}^{(1)}T_0(x_{j}) \\
     						 &=w_{00}^{(1)}{T_0(x_{0})}+w_{01}^{(1)}{T_0(x_{1})}+\ldots +w_{0N}^{(1)}{T_0(x_{N})} \\
								 &=w_{00}^{(1)}\delta_{00}+w_{01}^{(1)}\delta_{01}+\ldots +w_{0N}^{(1)}\delta_{0N} \\
							  0&=w_{00}^{(1)} 		 
\end{aligned}
\label{a00}
\end{equation*}

\noindent
Similarly, substituting the next basis function $T_1(x)$ into the differential quadrature approximation (\ref{dqm}) yields
\scriptsize
\begin{equation*}
\begin{aligned}
T_1^{\prime}(x_0)&=\sum \limits_{j=0}^{N}w_{0j}^{(1)}T_1(x_{j}) \\
     						              &=w_{00}^{(1)}{T_1(x_{0})}+w_{02}^{(1)}{T_1(x_{1})}+\ldots +w_{0N}^{(1)}{T_1(x_{N})} \\
								              &=w_{00}^{(1)}\delta_{10}+w_{01}^{(1)}\delta_{11}+\ldots +w_{0N}^{(1)}\delta_{1N} \\
\dfrac{\dfrac{\pi}{\Delta x}(x_0-1\Delta x)\cos{\dfrac{x_0-1\Delta x}{\Delta x} \pi}-\sin{\dfrac{x_0-1\Delta x}{\Delta x} \pi} }{\dfrac{\pi}{\Delta x}(x_0-1\Delta x)^2}&=w_{01}^{(1)} \\
\dfrac{\dfrac{\pi}{\Delta x}(0\Delta x-1\Delta x)\cos{\dfrac{0\Delta x-1\Delta x}{\Delta x} \pi}-\sin{\dfrac{0\Delta x-1\Delta x}{\Delta x} \pi} }{\dfrac{\pi}{\Delta x}(0\Delta x-1\Delta x)^2}&=w_{01}^{(1)}  \\
\dfrac{(0-1)\cos{(0-1)\pi}-\sin{(0-1)\pi} }{(0-1)^2(\Delta x)^2}&=w_{01}^{(1)} \\
\dfrac{\cos{(0-1)\pi}}{(0-1)\Delta x}&=w_{01}^{(1)} \\
\dfrac{(-1)^{(0-1)}}{(0-1)\Delta x}&=w_{01}^{(1)}
\end{aligned}
\label{a01}
\end{equation*}
\normalsize
where $\cos{(0-1)\pi}$ is replaced by $(-1)^{0-1}$. Substituting the basis function $T_m(x)$ into (\ref{dqm}) leads
\scriptsize
\begin{equation*}
\begin{aligned}
T_2^{\prime}(x_0)&=\sum \limits_{j=0}^{N}w_{0j}^{(1)}T_2(x_{j}) \\
     						              &=w_{00}^{(1)}{T_1(x_{0})}+w_{02}^{(1)}{T_1(x_{1})}+\ldots +w_{0N}^{(1)}{T_1(x_{N})} \\
								              &=w_{00}^{(1)}\delta_{20}+w_{01}^{(1)}\delta_{21}++w_{02}^{(1)}\delta_{22}\ldots +w_{0N}^{(1)}\delta_{1N} \\
\dfrac{\dfrac{\pi}{\Delta x}(x_0-2\Delta x)\cos{\dfrac{x_0-2\Delta x}{\Delta x} \pi}-\sin{\dfrac{x_0-2\Delta x}{\Delta x} \pi} }{\dfrac{\pi}{\Delta x}(x_0-2\Delta x)^2}&=w_{02}^{(1)} \\
\dfrac{\dfrac{\pi}{\Delta x}(0\Delta x-2\Delta x)\cos{\dfrac{0\Delta x-2\Delta x}{\Delta x} \pi}-\sin{\dfrac{0\Delta x-2\Delta x}{\Delta x} \pi} }{\dfrac{\pi}{\Delta x}(0\Delta x-2\Delta x)^2}&=w_{02}^{(1)}  \\
\dfrac{(0-2)\cos{(0-2)\pi}-\sin{(0-2)\pi} }{(0-2)^2(\Delta x)^2}&=w_{02}^{(1)} \\
\dfrac{\cos{(0-2)\pi}}{(0-2)\Delta x}&=w_{02}^{(1)} \\
\dfrac{(-1)^{(0-2)}}{(0-2)\Delta x}&=w_{02}^{(1)}
\end{aligned}
\label{a02}
\end{equation*} 
\normalsize
to give the related weight in an explicit form. Following the same methodology leads an explicit formulation for an arbitrary weight $w_{mj}^{(1)}$ as
\begin{equation*}
\begin{aligned}
w_{mj}^{(1)}&=\dfrac{(-1)^{m-j}}{\Delta x(m-j)},m\neq j \\
w_{mm}^{(1)}&=0  
\end{aligned}\label{amj}
\end{equation*}%
for the weights $w_{mj}^{(1)}$ that are related to the grid $x_m$.
\section{Discretization, Implementation of Boundary Conditions and Time integration}

\noindent
Approximating the derivative terms in the coupled system of partial differential equations (\ref{nlss}) gives the ordinary system of equations 
\begin{equation}
\begin{aligned}
\left. \dfrac{\partial g(x,t)}{\partial t}\right|_{x=x_m} &=\sum\limits_{j=0}^{N}{w_{m,j}^{(2)}f(x_j,t)}+\kappa\left(f^2(x_m,t)+g^2(x_m,t)\right)f(x_m,t) \\
\left. \dfrac{\partial f(x,t)}{\partial t}\right|_{x=x_m} &=-\sum\limits_{j=0}^{N}{w_{m,j}^{(2)}g(x_j,t)}-\kappa\left(f^2(x_m,t)+g^2(x_m,t)\right)g(x_m,t) \\
m&=0,1,\ldots , N
\end{aligned} \label{dnlss}
\end{equation}
related to each grid point $x_m$. Implementation of the homogeneous Dirichlet boundary conditions (\ref{bc}) converts the system (\ref{dnlss}) to 
\begin{equation}
\begin{aligned}
\left. \dfrac{\partial g(x,t)}{\partial t}\right|_{x=x_m} &=\sum\limits_{j=1}^{N-1}{w_{m,j}^{(2)}f(x_j,t)}+\kappa\left(f^2(x_m,t)+g^2(x_m,t)\right)f(x_m,t) \\
\left. \dfrac{\partial f(x,t)}{\partial t}\right|_{x=x_m} &=-\sum\limits_{j=1}^{N-1}{w_{m,j}^{(2)}g(x_j,t)}-\kappa\left(f^2(x_m,t)+g^2(x_m,t)\right)g(x_m,t) \\
m&=1,2, \ldots , N-1
\end{aligned} \label{dnlss1}
\end{equation}

\noindent
In order to integrate the space discretized system (\ref{dnlss1}), a bunch of methods of various orders covering Heun's method(HEUN), the classical Runge-Kutta methods of order from two to four (RK2,RK3,RK4) and some variations of higher order Runge-Kutta methods such as the Runge-Kutta-Fehlberg (RKF) and the Cash-Karp (CK) methods are used. The solutions are computed for different choices of discretization parameters $\Delta x$ and $\Delta t$ without any linearization. 

\section{On the Matrix Stability of the Proposed Methods}

\begin{figure}[hp]
\centering
    \subfigure[Order $2$]{
   \includegraphics[scale =0.65] {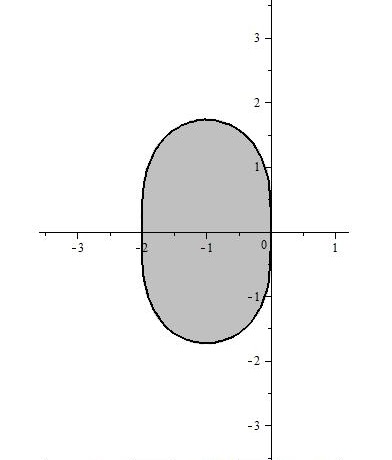}
   \label{fig:1a}
 }
   \subfigure[Order $3$]{
   \includegraphics[scale =0.65] {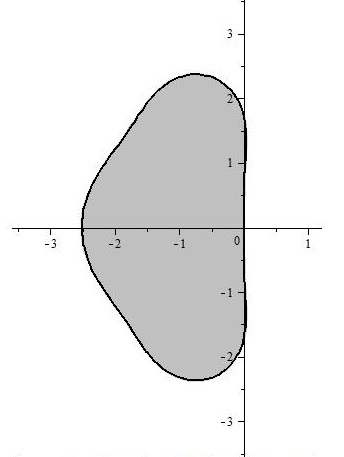}
   \label{fig:1b}
 }
 \subfigure[Order $4$]{
   \includegraphics[scale =0.65] {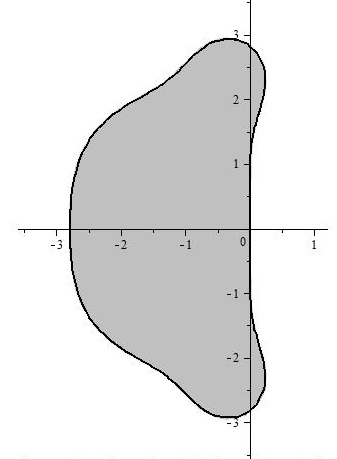}
   \label{fig:1c}
 }
 \subfigure[Order $5$]{
   \includegraphics[scale =0.65] {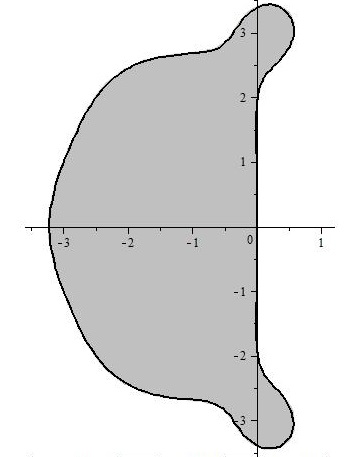}
   \label{fig:1d}
 }
 \caption{Stability regions for Runge-Kutta methods in the complex plane}
\end{figure}

\noindent
In this part, some basic concepts related to the stability of explicit Runge-Kutta methods are presented. The stability regions for each method used in this study are graphed in order to indicate the relation between the eigenvalue distribution and the stability of the method. The spectrum of the coefficient matrix gives the idea to choose appropriate $\Delta t$ to provide the stability condition. The stability region of a Runge-Kutta method of order $p$ depends on the stability polynomial $|S(\lambda \Delta t)|<1$ defined as
\begin{equation}
S(\lambda \Delta t) = \sum\limits_{s=0}^{p}{\dfrac{(\lambda \Delta t)^s}{s!}}
\end{equation}
since the common approximation in the Runge-Kutta methods are of the form
\begin{equation}
y_{n+1}=S(\lambda \Delta t) y_n
\end{equation}
for the differential equation
\begin{equation}
\dfrac{dy}{dt}=\lambda y
\end{equation}
The stability regions of the Runge-Kutta methods of orders from two to five are plotted in Fig \ref{fig:1a}-Fig \ref{fig:1d}.

\noindent
Now consider the linear equation system of order one such that 
\begin{equation}
\mathbf{\dfrac{dy}{dt}}=\mathbf{By}
\end{equation}
where $\mathbf{B}$ is the coefficient matrix with constant coefficients, $\mathbf{y}=\left[ y_1,y_2,\ldots y_N\right]^T$. The stability of this system is dependent on the eigenvalue distribution of the coefficient matrix $\mathbf{B}$.

\noindent
The space discretized system (\ref{dnlss1}) can be written in the matrix format as
\begin{equation}
\begin{bmatrix}g_1(t)\\ g_2(t) \\ \vdots \\ g_{N-1}(t)\\ f_1(t)\\f_2(t)\\ \vdots \\ f_{N-1}(t)\end{bmatrix}_t=\begin{bmatrix} \mathbf{0}& \mathbf{A}  \\ \mathbf{-A} & \mathbf{0} \end{bmatrix}\begin{bmatrix}g_1(t)\\ g_2(t) \\ \vdots \\ g_{N-1}(t)\\ f_1(t)\\f_2(t)\\ \vdots \\ f_{N-1}(t)\end{bmatrix} \label{ev}
\end{equation}
, 
\begin{equation}
\mathbf{A}=\begin{bmatrix}w_{11}^{(2)}+\kappa \tilde{\kappa_1} & w_{12}^{(2)} & \ldots & w_{1N-1}^{(2)} \\
w_{21}^{(2)}& w_{22}^{(2)}+\kappa \tilde{\kappa_2}  & \ldots & w_{2N-1}^{(2)} \\
\vdots & \ddots &&\vdots \\
w_{N-11}^{(2)} & w_{N-12}^{(2)} & \ldots & w_{N-1N-1}^{(2)}++\kappa \tilde{\kappa_{N-1}} 
\end{bmatrix}
\end{equation}
where $\tilde{\kappa_m}=f_m^2+g_m^2$ is assumed locally constant and $f_m$ and $g_m$ stand for $f(x_m,t)$ and $g(x_m,t)$, respectively. We investigate the eigenvalues of the coefficient matrix (\ref{ev}) by using the initial values of $f_m$ and $g_m$ for each test problem given in the following.
\section{Numerical Examples}
In this section of the study, some analytical and non-analytical initial-boundary value problems are considered. The error between the numerical and the analytical solutions is determined by using the discrete maximum error norm at the time $t$ defined as
\begin{equation*}
L_{\infty}(t)=\max\limits_{m}^{}\left|U_m-u_m\right|
\end{equation*}
where $u_m$ and $U_m$ are the analytical and numerical solutions at $x=x_m$, respectively.

\noindent
The lowest two conservation laws for the NLS equation defined as
\begin{equation}
\begin{aligned}
C_1&=\int\limits_{-\infty}^{\infty}{\left|u \right|^2 dx} \\
C_3&=\int\limits_{-\infty}^{\infty}{\left[\left|u_x \right|^2-\dfrac{1}{2}\kappa \left|u\right|^4\right] dx}
\end{aligned}
\end{equation}
are expected to remain constant as time proceeds\cite{zakharov1}. Even though $C_2$ and $C_4$ are also defined in the same paper, we do not calculate them here due to having complex components inside the integrals. Reporting the absolute relative change of the conservation laws at a specific time $t$ defined as 
\begin{equation}
\begin{aligned}
C(C_{\eta}(t))&=\left | \frac{C_{\eta}(t)-C_{\eta}(0)}{C_{\eta}(0)} \right | , \eta=1,3\\
\end{aligned}
\end{equation}
where $C_{\eta}(0),\eta=1,3$ denote the determined conservation law values initially can more useful to observe the preservation of the conservation laws. 
\subsection{Propagation of a Single Soliton}
\begin{figure}[H]
    \includegraphics[scale =0.5] {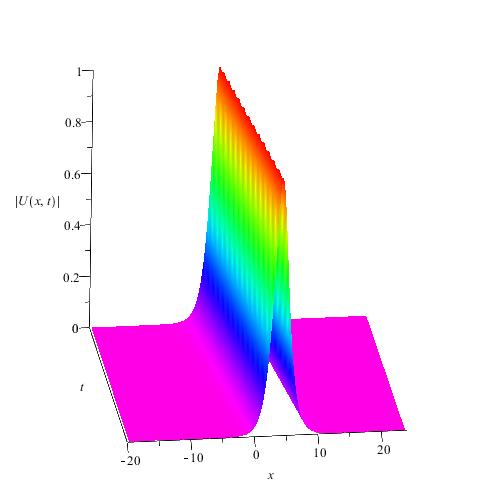}
 \caption{Propagation of single soliton simulation achieved by using parameters $\Delta t=0.0001$ and $\Delta x=0.1$}\label{SS}
\end{figure}

\noindent
The soliton solutions of the NLS is completely different from the solitons of KdV equation since they include both exponential and hyperbolic functions. These soliton solutions of the form
\begin{equation}
u(x,t)=\alpha\sqrt{\frac{2}{\kappa}}\exp{i\left[\frac{cx}{2}-\frac{(c^2-\alpha^2)t}{4}\right]}\sech{\alpha(x-ct)}
\end{equation}
comes into existence when the nonlinear and the dispersion terms are balanced completely. In this solution, $c$ stands for the velocity of the propagating soliton. In order to enable a comparison with some earlier studies, the parameters are selected as $c=4$, $\kappa =2$ and $\alpha=1$. Thus, the resultant envelop solution
\begin{equation}
\left|u(x,t)\right|=\sech{(x-4t)}
\end{equation}
represents a single soliton of constant speed $4$ propagating along the $x$-axis with height $1$. Various discretization parameters are used for the numerical illustrations in the experiment interval $[-20,24]$ up to the time $t=1$. The simulation of the propagation of the single soliton is depicted in Fig \ref{SS}.

\noindent
The time integration is carried out by using various techniques of orders from lower ones to higher ones. The discrete maximum error norms are tabulated in Table \ref{t1} for each choice of methods, space and time step sizes. 

\noindent
The initial values of conservation laws $C_1$ and $C_3$ are determined as $2$ and $7.33333333$ analytically. A comparison of discrete maximum errors and the absolute relative changes of the conservation laws are documented in Table \ref{t1}. 

\noindent
When $\Delta t=0.1$, all time integration techniques from HEUN to CK fail independent on the space grid sizes $\Delta x=0.3125$ and $\Delta x=0.1$ to generate the numerical solutions due to not satisfying the stability conditions. A discussion on the details of the stability is given in the following parts of this section. Since both HEUN and RK2 are of order two, the results generated by both methods are almost at the same digit-accuracy. When $\Delta x=0.3125$, the accuracies are determined in three, five and six decimal digits dependent on the size of time step $0.01$, $0.001$ and $0.0001$, respectively. Unfortunately, the decrease of the grid size to $0.1$ unbalances the stability for $\Delta t=0.1,0.01,0.001$. This balance is caught again with higher accuracy for $\Delta t=0.0001$. The accuracy reaches seven decimal digits with this choice of $\Delta t$. 

\noindent
The RK3 method gives five decimal-digit accuracy with $\Delta x=0.3125$ and $\Delta t=0.01$. Reducing $\Delta t$ to $0.001$ and $0.0001$ improves results to six decimal-digit accuracy. The decrease of $\Delta x$ to $0.1$ gives no results when $\Delta t=0.1,0.01$ but the choice of $\Delta t$ as $0.001$ and $0.0001$ gives eight and nine decimal-digit accuracies, respectively.

\noindent
The RK4, the RKF and the CK methods are accurate at six decimal digits with the discretization parameters $\Delta x=0.3125$ and $\Delta t=0.01,0.001,0.0001$. When the discretization parameters are chosen as $\Delta x=0.1$ and $\Delta t=0.001,0.0001$, nine decimal-digit accuracy are obtained in the results.

\noindent
The comparison with the earlier works indicates that the proposed methods solve this initial boundary value problem with acceptably higher accuracy. The proposed methods give better results than the results of the B-spline FEM\cite{dag1}, the B-spline Collocation \cite{gardner1} and the RBF methods\cite{dereli1}. Even though the results of CDQ\cite{korkmaz1} appear worse than other the present results in many cases the authors of that study explain the high error due to the forced boundary effect in $[-20,20]$. The change of the problem interval from $[-20,20]$ to $[-20,24]$ decreases the forced boundary effect on the accuracy of the results when Dirichlet conditions are used\cite{korkmaz1}. The results of the PDQ\cite{korkmaz2} are acceptably accurate but not better than the present results in many cases.

\noindent
The conservation laws calculated by using each numerical method are indicators of highly accurate results. The relative absolute changes decreases depending on the accuracy and the order of the method. Particularly, the absolute relative changes decreases below $10^{-10}$ when the stable higher order methods are used. The observations are in a good agreement with the theoretical expectations and the results reported in \cite{korkmaz1,korkmaz2}.

\begin{table}[hp]
\scriptsize
\caption{Comparison of the results with the results of some earlier works}
\begin{tabular}{lrrrrr}
\hline\hline
   Method  &  $\Delta x$ &  $\Delta t$ & $L_{\infty}(1)$ &  $C(C_1(1))$ & $C(C_3(1))$ \\
		\hline \hline
  HEUN (present)& $0.3125$ &  $0.1$ & $\infty$ &  $\infty$ & $\infty$ \\
       &  & $0.01$ &$2.629\times 10^{-3}$&$1.527\times 10^{-4}$&$9.563\times 10^{-4}$ \\
       && $0.001$ &$2.494\times 10^{-5}$&$1.215\times 10^{-7}$&$2.360\times 10^{-7}$ \\
       & & $0.0001$ &$1.378\times 10^{-6}$&$5.000\times 10^{-10}$&$1.363 \times 10^{-9}$\\

       &$0.1$ & $0.1$ &$\infty$ &$\infty$ &$\infty$  \\
       & & $0.01$ &$\infty$ &$\infty$ &$\infty$  \\
       & & $0.001$ &$\infty$ &$\infty$ &$\infty$  \\
       & & $0.0001$ &$2.471\times 10^{-7}$ &$<10^{-10}$& $2.727\times 10^{-10}$ \\

   RK2 (present) & $0.3125$ &  $0.1$ & $\infty$ &  $\infty$ & $\infty$ \\
       &  & $0.01$ &$2.629\times 10^{-3}$&$1.527\times 10^{-4}$&$9.563\times 10^{-4}$ \\
       && $0.001$ &$2.494\times 10^{-5}$&$1.215\times 10^{-7}$&$2.360\times 10^{-7}$ \\
       & & $0.0001$ &$4.837\times 10^{-6}$&$5.000\times 10^{-10}$&$1.363 \times 10^{-9}$\\
       
       &$0.1$ & $0.1$ &$\infty$ &$\infty$ &$\infty$  \\
       & & $0.01$ &$\infty$ &$\infty$ &$\infty$  \\
       & & $0.001$ &$\infty$ &$\infty$ &$\infty$  \\
       & & $0.0001$ &$2.471\times 10^{-7}$ &$<10^{-10}$& $2.727\times 10^{-10}$ \\
            
   RK3 (present)& $0.3125$ &  $0.1$ & $\infty$ &  $\infty$ & $\infty$ \\
       &  & $0.01$ &$7.560\times 10^{-5}$&$4.064\times 10^{-5}$& $7.883\times 10^{-5}$ \\
       && $0.001$ &$1.378\times 10^{-6}$&$4.050\times 10^{-8}$&$7.854\times 10^{-8}$ \\
       & & $0.0001$ &$1.378\times 10^{-6}$&$5.000\times 10^{-10}$&$4.090\times 10^{-10}$ \\   
       &$0.1$ & $0.1$ &$\infty$ &$\infty$ &$\infty$  \\
       & & $0.01$ &$\infty$ &$\infty$ &$\infty$  \\
       & & $0.001$ &$7.510\times 10^{-8}$ &$4.100\times 10^{-8}$&$7.854\times 10^{-8}$  \\
      & & $0.0001$ &$2.800\times 10^{-9}$& $5.000\times 10^{-10}$ & $8.181\times 10^{-10}$  \\
   RK4 (present)& $0.3125$ &  $0.1$ & $\infty$ &  $\infty$ & $\infty$ \\
       &  & $0.01$ &$2.092\times 10^{-6}$&$4.500\times 10^{-8}$& $1.262\times 10^{-7}$\\
       && $0.001$ &$1.378\times 10^{-6}$&$5.000\times 10^{-10}$&$1.227\times 10^{-9}$ \\
       & & $0.0001$&$1.378\times 10^{-6}$&$5.000\times 10^{-10}$&$1.227\times 10^{-9}$ \\
  
       &$0.1$ & $0.1$ &$\infty$ &$\infty$ &$\infty$  \\
       & & $0.01$ &$\infty$ &$\infty$ &$\infty$  \\      
       & & $0.001$ &$2.814\times 10^{-9}$ &$<10^{-10}$ &$<10^{-10}$  \\
       & & $0.0001$ &$2.805\times 10^{-9}$ &$<10^{-10}$  & $<10^{-10}$ \\
       
 RKF (present)& $0.3125$ &  $0.1$ & $\infty$ &  $\infty$ & $\infty$ \\
       &  & $0.01$ &$1.380\times 10^{-6}$&$6.000 \times 10^{-9}$&$4.750 \times 10^{-8}$ \\
       && $0.001$ &$1.378\times 10^{-6}$&$< 10^{-10}$&$<10^{-10}$ \\
       & & $0.0001$ &$1.378\times 10^{-6}$&$<10^{-10}$&$<10^{-10}$ \\
       &$0.1$ & $0.1$ &$\infty$ &$\infty$ &$\infty$  \\
       & & $0.01$ &$\infty$ &$\infty$ &$\infty$  \\
       & & $0.001$ &$2.803 \times 10^{-9}$ &$<10^{-10}$&$<10^{-10}$  \\
       & & $0.0001$ &$2.805\times 10^{-9}$ &$<10^{-10}$ &$<10^{-10}$  \\
 
 CK (present)& $0.3125$ &  $0.1$ & $\infty$ &  $\infty$ & $\infty$ \\
       &  & $0.01$ &$1.379\times 10^{-6}$&$1.190\times 10^{-9}$& $1.808\times 10^{-9}$\\
       && $0.001$ &$1.378\times 10^{-6}$&$<10^{-10}$&$<10^{-10}$ \\
       & & $0.0001$ &$1.378\times 10^{-6}$&$<10^{-10}$&$<10^{-10}$ \\
       &$0.1$ & $0.1$ &$\infty$ &$\infty$ &$\infty$  \\
       & & $0.01$ &$\infty$ &$\infty$ &$\infty$  \\
       & & $0.001$ &$2.804\times 10^{-9}$ &$<10^{-10}$ &$<10^{-10}$  \\
       & & $0.0001$ &$2.805\times 10^{-9}$ &$<10^{-10}$ &$<10^{-10}$  \\
CDQ\cite{korkmaz1} & $0.3125$ &$0.025$ &$5.562\times 10^{-5}$ &$4.380\times 10^{-6}$&$1.214\times 10^{-5}$ \\
                  &          &$0.01$   &$2.089\times 10^{-6}$ &$4.512\times 10^{-8}$&$1.253\times 10^{-7}$ \\
                  &  $0.125$ &$0.0025$ &$1.557\times 10^{-6}$ &$4.412\times 10^{-11}$&$1.226\times 10^{-10}$ \\
                  &          &$0.001$  &$1.550\times 10^{-7}$ &$4.541\times 10^{-13}$&$1.266\times 10^{-12}$ \\
PDQ\cite{korkmaz2}& $0.3125$ &$0.02$ &$2.531\times 10^{-5}$ &$1.440\times 10^{-6}$&$3.941\times 10^{-6}$ \\
                  & $0.1$    &$0.0025$   &$1.932\times 10^{-7}$ &$5.222\times 10^{-11}$&$2.930\times 10^{-9}$ \\
B-spline FEM \cite{dag1}& $0.3125$ &$0.02$ &$0.002$ &$6.600\times 10^{-6}$&$3.417\times 10^{-4}$ \\
                        & $0.05$ &$0.005$ &$3.000\times10^{-4}$ &$<10^{-8}$&$6.000\times 10^{-7}$ \\
B-spline Collocation \cite{gardner1}&$0.05$ &$0.005$ &$0.008$ &$<10^{-6}$&$<10^{-6}$ \\ 
                                    &$0.03$ &$0.005$ &$0.002$ &$<10^{-6}$&$<10^{-6}$ \\
RBF Method G \cite{dereli1}&$0.3125$ &$0.001$ &$2.800\times 10^{-5}$ &&\\
RBF Method MQ \cite{dereli1}&$0.3125$ &$0.001$ &$2.165\times 10^{-3}$ &&\\
RBF Method IMQ \cite{dereli1}&$0.3125$ &$0.001$ &$4.860\times 10^{-4}$ &&\\                                     RBF Method IQ \cite{dereli1}&$0.3125$ &$0.001$ &$5.652\times 10^{-3}$ &&\\                                     
\hline\hline
\end{tabular}
\label{t1}
\end{table}

\begin{figure}[H]
 
   \subfigure[$\Delta x=0.3125$]{
   \includegraphics[scale =0.4] {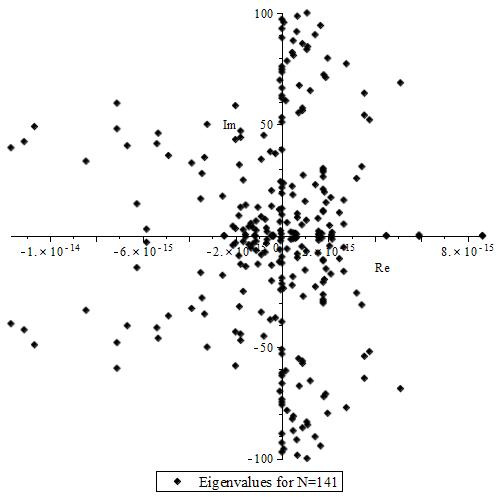}
   \label{ev13}
 }
 \subfigure[$\Delta x=0.3125$]{
   \includegraphics[scale =0.4] {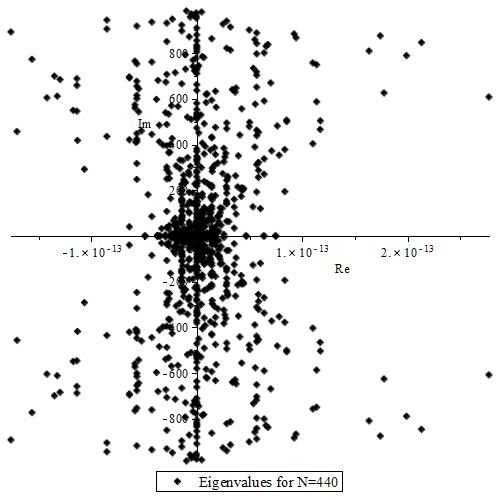}
   \label{ev11}
 }
 \caption{The eigenvalue distributions for propagation of a single soliton}
\end{figure}

\noindent
The stability analysis of the proposed methods depends on the eigenvalue distribution of the coefficient matrix given in the previous section. Assuming the real and the complex components of the solution are locally constant converts the system to a linear system with constant coefficient matrix $\mathbf{B}$. The spectrum of this coefficient matrix gives the information to control the size of $\Delta t$ for the stability of the method. When $\Delta x$ is chosen as $0.3125$, the complex components of the eigenvalues are larger than $100$ in absolute value, Fig \ref{ev13}. The choice of the time step size $\Delta t$ should squeeze the value of the multiplication of each eigenvalue $\lambda_i$ and $\Delta t$ to the regions given in Fig \ref{fig:1a}-Fig \ref{fig:1d} correspond to the order. Similar conclusion can be written for the choice of $\Delta x=0.1$ by commenting on Fig \ref{ev11}. Since some of the complex components of the eigenvalues are greater than $900$ in absolute value, the stability forces $\Delta t$ to be chosen smaller. This perspective explains why time integration techniques fail to solve the problem.

\subsection{Collision of Two Positive Solitary Waves Moving in the Opposite Directions}
\begin{figure}[H]
    \includegraphics[scale =0.5] {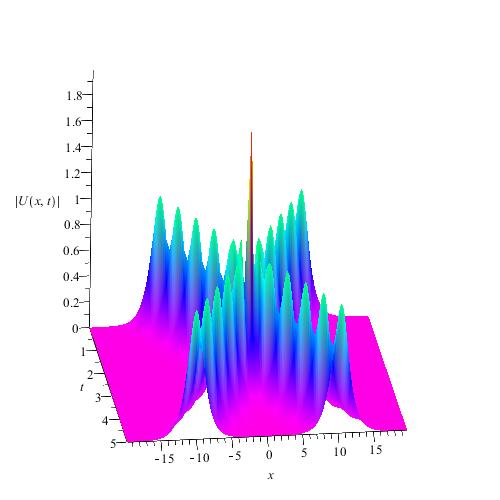}
 \caption{Collision of two positive solitary waves propagating in the opposite directions}\label{inter}
\end{figure}

The collision of two positive solitary waves moving in the opposite directions along the $x$-axis is modeled by the initial data 
\begin{equation*}
u(x,0)=\sum\limits_{j=1}^{2}{u_{j}(x,0)}
\end{equation*}
where
\begin{equation*}
u_j(x,0)=\beta_j\sqrt{\frac{2}{\kappa}}\exp{\left(i\frac{1}{2}c_j(x-\hat{x_j})\right)} \sech{\beta_j\left(x-\hat{x_j} \right)},\, j=1,2
\end{equation*}
where $\hat{x_j}$ stands for the peak positions of the solitary waves initially. The experiment is completed with the appropriate parameters $\kappa =2$, $\beta_j =1, j=1,2$, $c_1=4$ and $c_2=-4$ used in the previous studies \cite{korkmaz1,korkmaz2,dag1,gardner1}. The peaks of the well separated two positive solitary waves are sited to $x=-10$ and $x=10$ initially by choosing $\hat{x_1}=-10$ and $\hat{x_2}=10$, respectively. The numerical solutions are computed with the discretization parameters $\Delta x=0.25$ and $\Delta t=0.005$ in the artificial domain $[-20,20]$, Fig \ref{inter}. All routines are run up to time $t=5$ to observe the separation clearly after the collision. This choice of simulation ending time restricts both waves to hit the ends of artificial interval, too. The homogeneous Dirichlet boundary conditions are used to provide the compatibility with theory. 

\noindent
Both solitary waves of unit height initially move towards each other due to the difference in the sign of $c_1$ and $c_2$. The collision starts as time proceeds and can be observable clearly around the time $t=2$. When the time is $t=2.5$, the height of the joint waves exceeds $1.9$ because the structure of the solution is sum of two positive waves. Both waves keep moving along their own ways and separate from each other as time proceeds. After full separation, both waves turn their initial shapes and heights. Each one takes the initial position of the other one owing to their constant velocities $4$ when the time reaches $t=5$.

\begin{figure}[H]
    \includegraphics[scale =0.5] {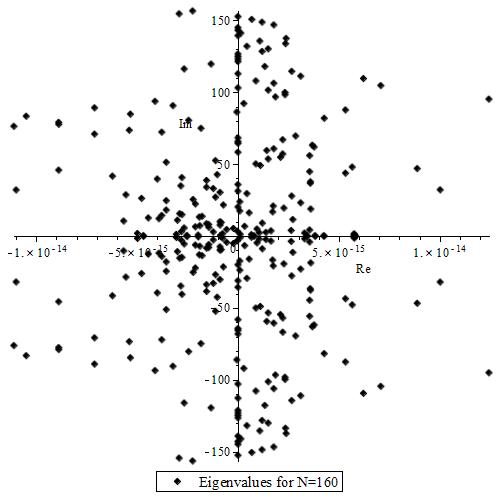}
 \caption{Eigenvalues for the collision of two positive solitaries problem}\label{interev}
\end{figure}

\noindent
The HEUN and the RK2 both fail to simulate the solutions with those discretization parameters. The analytical values of the conservation laws are $C_1=4$ and $C_3=14.66666667$ initially. The absolute relative changes of the conservation laws are the indicators of the accuracy of the method for this problem due to non existence of the analytical solution. The reported results of absolute relative changes indicate that $C_1$ and $C_3$ both change in five decimal digits when the method is RK3, Table \ref{t2}. The RK4 preserves the conservation laws $C_1$ and $C_3$ in nine and eight decimal digits, respectively. The first conservation law is changed relatively in absolute value in ten decimal digits for the methods RKF and CK. The performance of RKF is one decimal digit better than eight decimal digits preservation of CK.
\begin{table}[H]
\caption{Absolute relative changes of the conservation laws at $t=5$}
\begin{tabular}{lrrrr}
\hline\hline
   Method  &  $\Delta x$ &  $\Delta t$  &  $C(C_1(1))$ & $C(C_3(1))$ \\
		\hline \hline      
   RK3 (present)& $0.25$ &  $0.005$  &  $2.494\times 10^{-5}$ & $4.793\times 10^{-5}$ \\
   RK4 (present)& $0.25$ &  $0.005$  &  $6.250\times 10^{-9}$& $3.068\times 10^{-8}$ \\
   RKF (present)& $0.25$ &  $0.005$ &  $5.000\times 10^{-10}$ & $8.863\times 10^{-9}$ \\
   CK (present)& $0.25$ &  $0.005$ &  $5.000\times 10^{-10}$ & $1.090\times 10^{-8}$ \\
 PDQ\cite{korkmaz2} & $0.25$ &  $0.01$ &  $2.215\times 10^{-7}$ & $4.358\times 10^{-7}$ \\
\hline\hline
\end{tabular}
\label{t2}
\end{table}
\noindent
The eigenvalue distribution related to this problem indicates that the choice of $\Delta t=0.005$ is not sufficient to carry all $\lambda_j \Delta t$ given in Fig \ref{interev} to the stability region Fig \ref{fig:1a} for the second order methods, HEUN and RK2. The maximum and minimum complex components of all eigenvalues are $\pm 156.462$ and $\pm156.462 \times 0.005=\pm 0.782$ is at the outside of the stability region \ref{fig:1a}. However, that choice of $\Delta t$ is sufficient for the stability of the RK3, the RK4, the RKF and the CK methods. 

\subsection{Birth of Standing Initial Pulse}
\begin{figure}[H]
 
   \subfigure[$A=1$]{
   \includegraphics[scale =0.4] {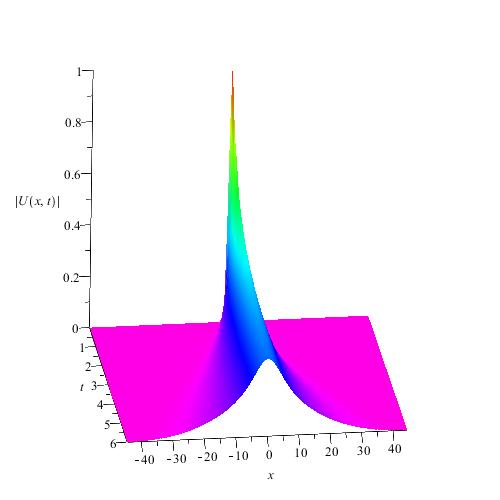}
   \label{BOSA1}
 }
 \subfigure[$A=1.78$]{
   \includegraphics[scale =0.4] {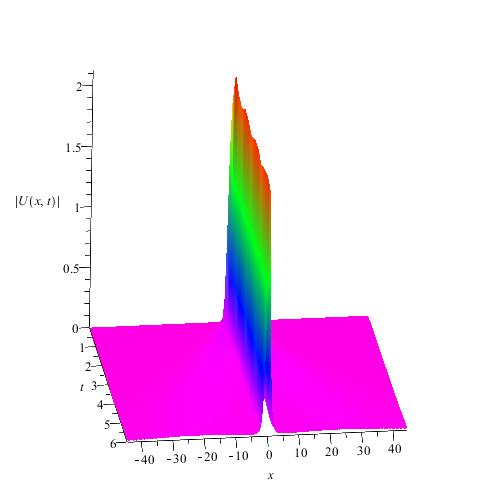}
   \label{BOSA178}
 }
 \caption{Fade out and birth of a standing initial pulse}
\end{figure}

\noindent
The birth of a standing pulse occurs when the integral of initial data is greater than $\pi$ in the infinite interval. Otherwise the initial pulse fades out as time proceeds. The Maxwellian initial condition \cite{gardner1}
\begin{equation}
u(x,0)=Ae^{-x^2}
\end{equation}
is chosen to demonstrate the experiment in the interval $[-45,45]$. The algorithms are run up to the ending time $t=6$ with $\kappa =2$, fixed grid size $\Delta x=0.5$ and various $\Delta t$ values.
The choice $A=1$ gives a positive pulse of unit height positioned at $x=0$ initially but this value of $A$ does not keep the balance in the equation and the initial pulse fades out as time proceeds, Fig \ref{BOSA1}. On the other hand, the choice of $A$ as $1.78$ gives an initial pulse of height $1.78$. At the earlier times of the motion, the height increases rapidly and gets larger than $2$ as giving birth from the both sides of its basis, Fig \ref{BOSA178}. As these two small pulses goes far away from the initial pulse, the height of their mother decreases to below $2$. The mother keeps its position, shape and height as the births propagates.

\begin{figure}[H]
 
   \subfigure[$A=1$]{
   \includegraphics[scale =0.4] {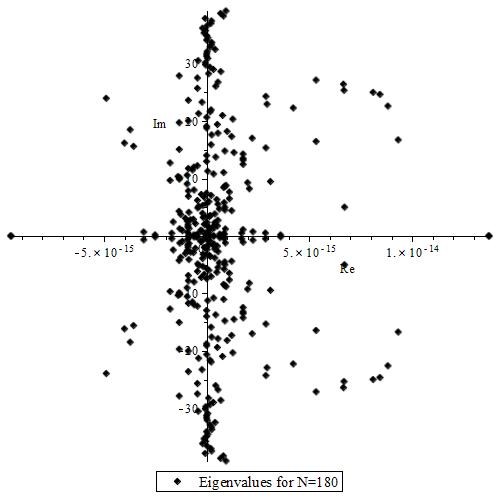}
   \label{EVA1}
 }
 \subfigure[$A=1.78$]{
   \includegraphics[scale =0.4] {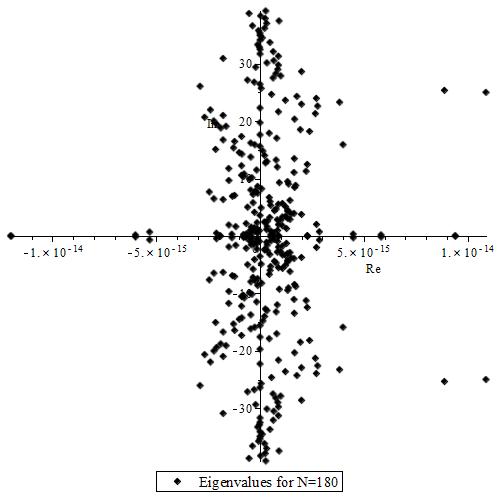}
   \label{EVA178}
 }
 \caption{Eigenvalue  distribution for the birth of standing initial pulse}
\end{figure}

\noindent
The conservation laws are computed from the initial data as $C_1=A^2\sqrt{\pi/2}\cong 3.97100$ and $C_3=(\sqrt{\pi}/4)A^2(2\sqrt{2}-\kappa A^2) \cong -4.92653$. The absolute relative changes of the conservation laws are tabulated in Table \ref{t3}. The discretization parameter choice $\Delta x=0.5$ and $\Delta t=0.05$ causes instability for the second and third order methods. Even though all other choices of the discretization parameters generate acceptably well solutions, the absolute relative change of the third conservation law does not improve for smaller time step sizes. On the contrary, the first conservation law is sensitive to smaller time step sizes when $\Delta x$ is fixed as $0.5$ and improves in digits. The accuracies do not improve too much even the order of the method increases. The comparison with some earlier differential quadrature methods indicates that the present results are satisfactorily well. 
\begin{table}[H]
\caption{Absolute relative changes of the conservation laws at $t=6$ for the birth of standing pulse}
\begin{tabular}{lrrrr}
\hline\hline
   Method  &  $\Delta x$ &  $\Delta t$  &  $C(C_1(6))$ & $C(C_3(6))$ \\
		\hline \hline
		HEUN (present)& $0.5$ &  $0.05$  &  $\infty$ & $\infty$ \\
		HEUN (present)& $0.5$ &  $0.005$  &  $6.697\times 10^{-5}$ & $7.708\times 10^{-4}$ \\
		HEUN (present)& $0.5$ &  $0.0005$  &  $6.623\times 10^{-8}$ & $7.362\times 10^{-5}$ \\
		RK2 (present)& $0.5$ &  $0.05$  &  $\infty$ & $\infty$ \\
		RK2 (present)&  $0.5$ &  $0.005$  &  $6.697\times 10^{-5}$ & $7.708\times 10^{-4}$ \\
		RK2 (present)&  $0.5$ &  $0.0005$  &  $6.623\times 10^{-8}$ & $7.362\times 10^{-5}$ \\
	 RK3 (present)& $0.5$ &  $0.05$  &  $\infty$ & $\infty$ \\
   RK3 (present)& $0.5$ &  $0.005$  &  $2.207\times 10^{-5}$ & $7.163\times 10^{-4}$ \\
   RK3 (present)& $0.5$ &  $0.0005$  &  $2.241\times 10^{-8}$ & $7.361\times 10^{-4}$ \\
   RK4 (present)& $0.5$ &  $0.05$  &  $1.267\times 10^{-4}$& $4.179\times 10^{-3}$ \\
   RK4 (present)& $0.5$ &  $0.005$  &  $2.770\times 10^{-9}$& $7.361\times 10^{-4}$ \\
   RK4 (present)& $0.5$ &  $0.0005$  &  $7.554\times 10^{-10}$& $7.361\times 10^{-4}$ \\
   
   RKF (present)& $0.5$ &  $0.05$  &  $1.571\times 10^{-4}$& $4.153\times 10^{-3}$ \\
   RKF (present)& $0.5$ &  $0.005$  &  $1.259\times 10^{-9}$& $7.361\times 10^{-4}$ \\
   RKF (present)& $0.5$ &  $0.0005$  &  $7.554\times 10^{-10}$& $7.361\times 10^{-4}$ \\
   
   CK (present)& $0.5$ &  $0.05$  &  $5.548\times 10^{-5}$& $8.082\times 10^{-4}$ \\
   CK (present)& $0.5$ &  $0.005$  &  $1.259\times 10^{-9}$& $7.361\times 10^{-4}$ \\
   CK (present)& $0.5$ &  $0.0005$  &  $7.554\times 10^{-10}$& $7.361\times 10^{-4}$ \\ 
   
CDQ\cite{korkmaz1} & $0.25$ &  $0.01$ &  $8.836\times 10^{-8}$ & $3.892\times 10^{-6}$ \\   
PDQ\cite{korkmaz2} & $0.25$ &  $0.01$ &  $6.334\times 10^{-7}$ & $7.110\times 10^{-5}$ \\

\hline\hline
\end{tabular}
\label{t3}
\end{table}
The eigenvalues for $A=1$ and $A=1.78$ are graphed in Fig \ref{EVA1} and Fig \ref{EVA178}, respectively. Since the maximum and minimum complex components of all eigenvalues are $\pm 39.155$, the choice of $\Delta t=0.05$ or greater values is not sufficient for the stability, Fig \ref{fig:1a} and Fig \ref{fig:1b}. The multiplication of the maximum( and the minimum) complex component is equal to $\pm 1.957$. These values are sufficient to satisfy the stability condition for only the higher order methods RK4, RKF and CK.  

\subsection{Birth of Propagating Pulse}
\begin{figure}[H]
 
   \subfigure[$A=1$]{
   \includegraphics[scale =0.4] {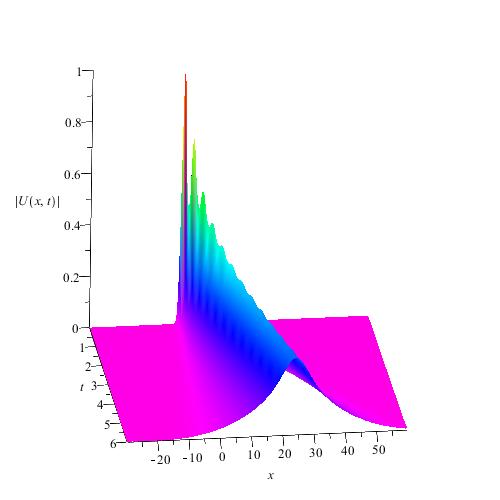}
   \label{BOP1}
 }
 \subfigure[$A=1.78$]{
   \includegraphics[scale =0.4] {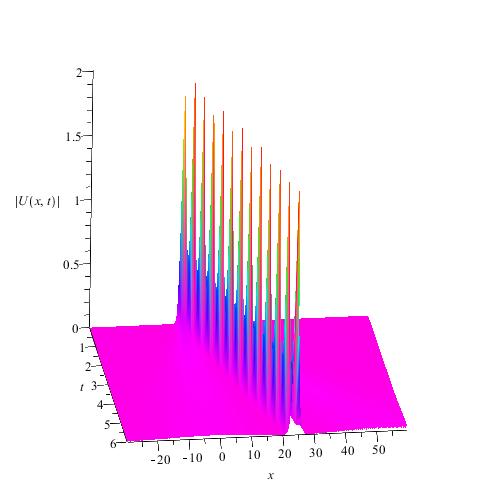}
   \label{BOP178}
 }
 \caption{The behaviors of the propagating pulse depending on $A$}
\end{figure}

\noindent
The birth of propagating pulse problem is a model with the initial data
\begin{equation}
u(x,0)=Ae^{-x^2+2ix}
\end{equation}
where the constant $A$ denotes the height of the initial pulse. This initial pulse propagates along the $x$-axis as time proceeds. The choice $A=1$ causes the initial pulse of unit height to fade out depending on the proceeding time, Fig \ref{BOP1}. The designed routines are run up to the time $t=6$ with a fixed $\Delta x=0.25$ and $\kappa=2$ over the problem interval $[-30,60]$ to simulate the solutions illustrating the birth of propagating pulse for $A=1.78$, Fig \ref{BOP178}. This choice of $A$ produces an initial pulse of height $1.78$. The height increases over $2$ and two bulges begin to appear at both sides of the pulse at the earlier times of the solution. These two bulges moves far away from the pulse. Subsequently, the height comes below $2$ and stays fixed around $1.9$. The two bulges moving far away from the pulse lose their heights but widen as time goes.  

\noindent
The conservation laws are calculated initially as $C_1=3.971000512$ and $C_3=10.95838443$. The second order methods HEUN and RK2 fail to give the solutions when $\Delta t=0.01$ owing to the instability, Table \ref{t4}. The reduction of $\Delta t$ to $0.001$ resolve this issue by giving six decimal-digit preservation in the $C_1$ and four decimal-digit preservation in the $C_3$ at the simulation ending time $t=6$ for both methods. The RK3 generates three decimal-digit preservation in both of the conservation laws with the discretization parameters $\Delta x=0.25$ and $\Delta t=0.01$. When the time step size is reduced to $0.001$, the preservation of the absolute relative changes of both $C_1$ and $C_3$ improve to six and four decimal digits, respectively. $\Delta t=0.01$ gives five, six and seven decimal-digit preservation in the values of $C_1$ for the methods RK4, RK5 and CK, again respectively. This choice of $\Delta t$ preserves the $C_3$ in four decimal digits for all the higher order methods RK4, RKF and CK. Even though decreasing the time step size to $0.001$ improve the preservation in $C_1$ more than ten decimal digits, no decimal digit improvement in the preservation is observed in the value of the $C_3$. 

\begin{figure}[H]
 
   \subfigure[$A=1$]{
   \includegraphics[scale =0.4] {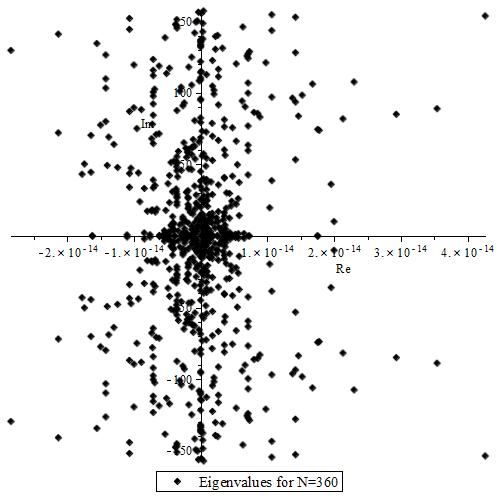}
   \label{EV41}
 }
 \subfigure[$A=1.78$]{
   \includegraphics[scale =0.4] {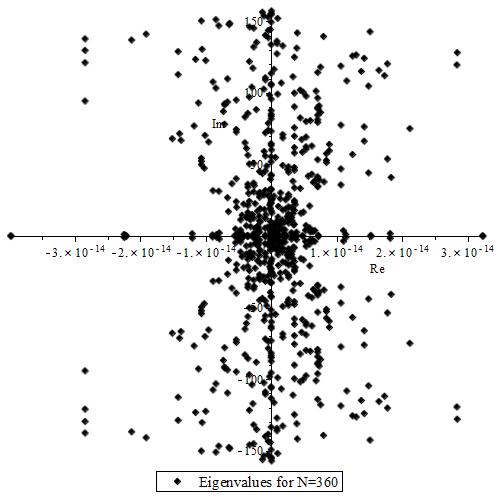}
   \label{EV4178}
 }
 \caption{The behaviors of the propagating pulse depending on $A$}
\end{figure}

\noindent
The eigenvalue distributions of coefficient matrices are graphed in Fig \ref{EV41} - \ref{EV4178}. The maximum complex component in absolute value of all eigenvalues is $157.267$ requires smaller $\Delta t$ values in the lower order methods HEUN and RK2. $\lambda_j \times \Delta t$ is not sufficient to satisfy the stability condition given in Fig \ref{fig:1a} when $\Delta t=0.01$. The choices of time step size for the other methods used in this study supply with the stability condition.
\begin{table}[H]
\caption{Absolute relative changes of the conservation laws at $t=6$ for the birth of propagating pulse}
\begin{tabular}{lrrrr}
\hline\hline
   Method  &  $\Delta x$ &  $\Delta t$  &  $C(C_1(6))$ & $C(C_3(6))$ \\
		\hline \hline
		HEUN (present)& $0.25$ &  $0.01$  &  $\infty$ & $\infty$ \\
		HEUN (present)& $0.25$ &  $0.001$  &  $4.474\times 10^{-6}$ & $1.835\times 10^{-4}$ \\
		
		RK2 (present)& $0.25$ &  $0.01$  &  $\infty$ & $\infty$ \\
		RK2 (present)& $0.25$ &  $0.001$  &  $4.474\times 10^{-6}$ & $1.835\times 10^{-4}$ \\	
	 
	 RK3 (present)& $0.25$ &  $0.01$  &  $1.314\times 10^{-3}$ & $4.865\times 10^{-3}$ \\
   RK3 (present)& $0.25$ &  $0.001$  &  $1.491\times 10^{-6}$ & $2.113\times 10^{-4}$ \\
  
   RK4 (present)& $0.25$ &  $0.01$  &  $1.462\times 10^{-5}$& $3.862\times 10^{-4}$ \\
   RK4 (present)& $0.25$ &  $0.001$  &  $< 10^{-10}$& $2.041\times 10^{-4}$ \\

   RKF (present)& $0.25$ &  $0.01$  &  $1.779\times 10^{-6}$& $1.816\times 10^{-4}$ \\
   RKF (present)& $0.25$ &  $0.001$  &  $< 10^{-10}$& $2.041\times 10^{-4}$ \\

   CK (present)& $0.25$ &  $0.01$  &  $2.727\times 10^{-7}$& $1.995\times 10^{-4}$ \\
   CK (present)& $0.25$ &  $0.001$  &  $< 10^{-10}$& $2.041\times 10^{-4}$ \\

CDQ\cite{korkmaz1} & $0.25$ &  $0.01$ &  $1.461\times 10^{-5}$ & $1.824\times 10^{-4}$ \\   
PDQ\cite{korkmaz2} & $0.25$ &  $0.01$ &  $1.575\times 10^{-5}$ & $5.265\times 10^{-4}$ \\

\hline\hline
\end{tabular}
\label{t4}
\end{table}

\subsection{Bound State of $n$ solitons}
The bound state of $n$ solitons for the nonlinear cubic NLS equation (\ref{nls}) is simulated with the initial data
\begin{equation}
u(x,0)=\sech{x}
\end{equation}
where $\kappa \geq 2n^2$, $n \in Z^{+}$ in the NLS equation (\ref{nls})\cite{zakharov1,miles1}. When this condition on $\kappa$ is satisfied, the $u(x,t)$ gives $n$-soliton bound state type solution and the periodicity property in a mobile frame of the solution is fulfilled. The case $2n^2<\kappa<2(n+1)^2$ causes the asymptotic solution to include both bound state and an oscillation with $O(t^{-1/2})$ as $\kappa<2$ gives only the oscillation\cite{miles1}. Unfortunately, the larger values of $n$ is not sufficiently suitable for the numerical solutions. In the present study, the behaviors are investigated for $n=2,3,4$ in the interval $[-20,20]$. The designed algorithms are run by using the discretization parameters $\Delta x=0.125$ and $\Delta t=0.001$ up to the ending time $t=0.6$. 

\begin{figure}[H]
\centering
   \subfigure[$t=0$]{
   \includegraphics[scale =0.3] {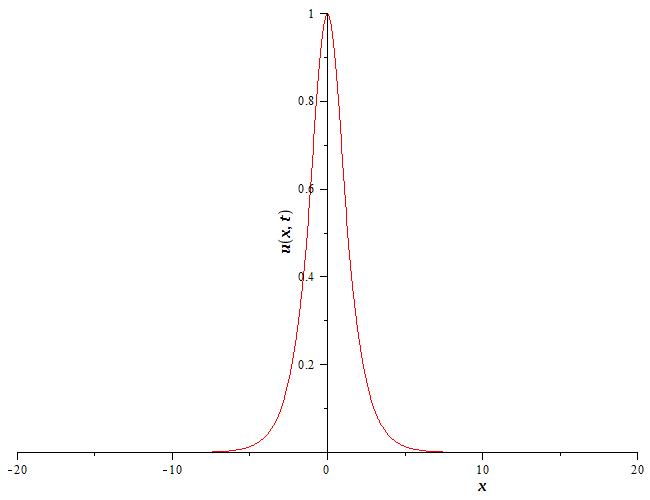}
   \label{bound20}
 }
   \subfigure[$t=0.2$]{
   \includegraphics[scale =0.3] {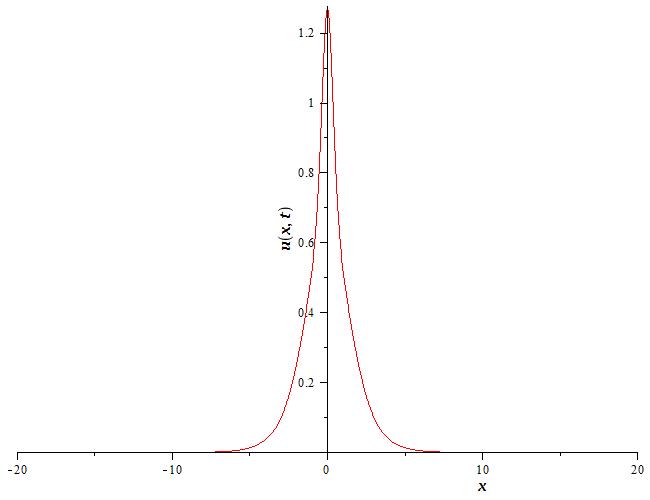}
   \label{bound22}
 }
    \subfigure[$t=0.3$]{
   \includegraphics[scale =0.3] {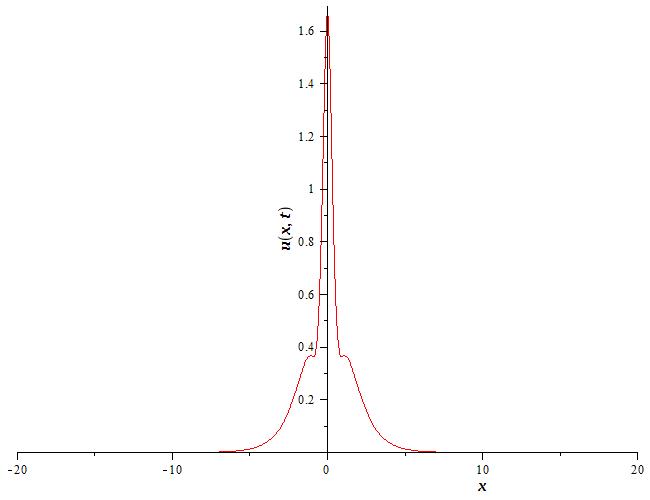}
   \label{bound23}
 }
    \subfigure[$t=0.4$]{
   \includegraphics[scale =0.3] {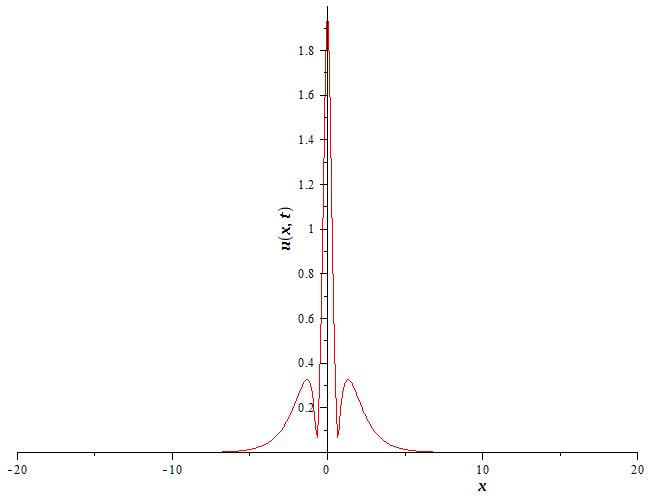}
   \label{bound24}
 }
    \subfigure[$t=0.5$]{
   \includegraphics[scale =0.3] {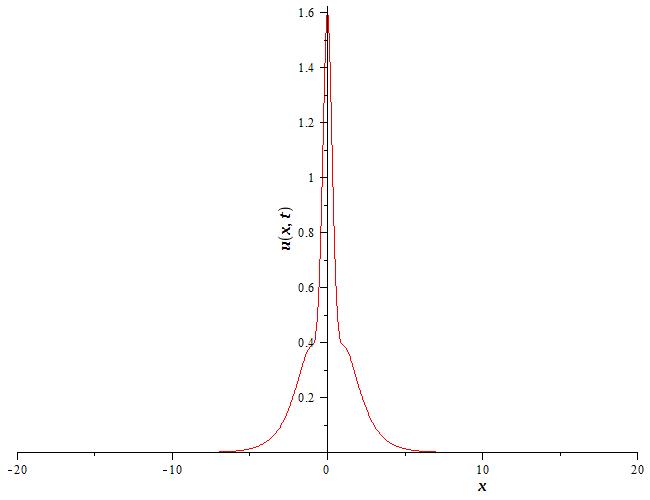}
   \label{bound25}
 }
    \subfigure[$t=0.6$]{
   \includegraphics[scale =0.3] {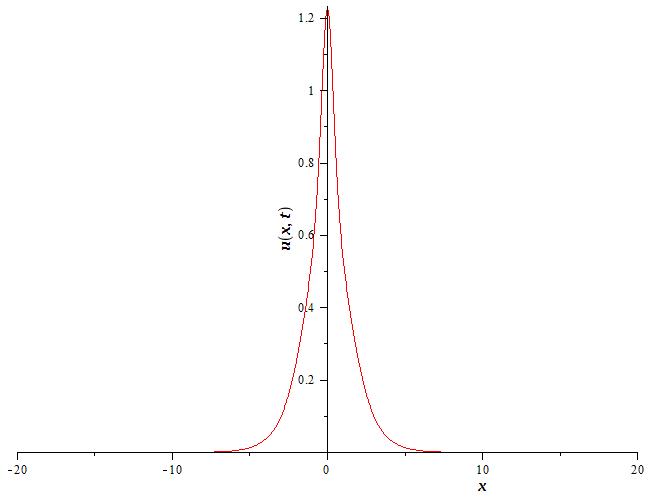}
   \label{bound26}
 }
\caption{Bound state of solitons for $n=2$}
\end{figure}

\noindent
In the first experiment, $n$ is chosen as $2$. The initial data stands for a pulse of unit height positioned at $x=0$, Fig \ref{bound20}. When the time reaches $t=0.2$ the height of the pulse passes over $1.2$, Fig \ref{bound22}. Two bulges begin to appear at both sides of this pulse when the time is $t=0.3$, Fig \ref{bound23}, as the peak gets higher than $1.6$. These two bulges are separated from the pulse at the time $t=0.4$, Fig \ref{bound24}. As the time proceeds, the bulges merges to the pulse again at $t=0.5$, Fig \ref{bound25}. However, the height of the pulse descends down $1.6$. At the end of the simulation, the height is measured just over $1.2$ and the bulges disappear, Fig \ref{bound26}. 

\begin{figure}[H]
\centering
   \subfigure[$t=0.175$]{
   \includegraphics[scale =0.3] {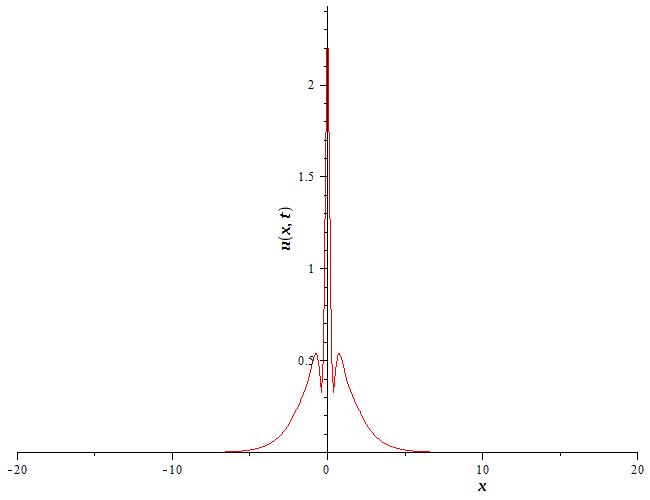}
   \label{bound38}
 }
   \subfigure[$t=0.250$]{
   \includegraphics[scale =0.3] {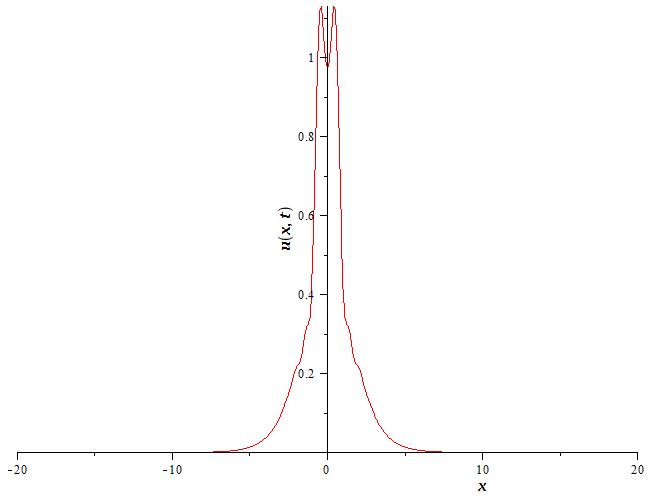}
   \label{bound311}
 }
    \subfigure[$t=0.300$]{
   \includegraphics[scale =0.3] {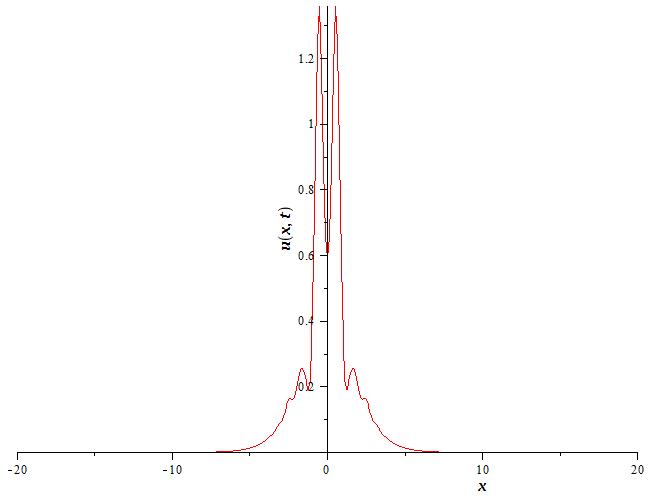}
   \label{bound313}
 }
    \subfigure[$t=0.325$]{
   \includegraphics[scale =0.3] {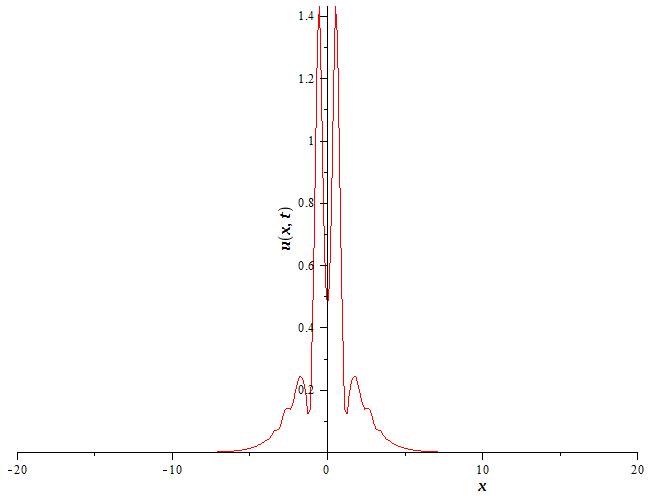}
   \label{bound314}
 }
    \subfigure[$t=0.400$]{
   \includegraphics[scale =0.3] {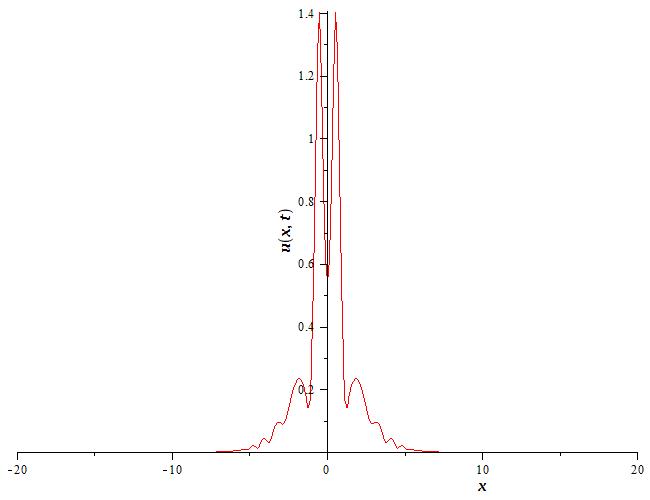}
   \label{bound317}
 }
    \subfigure[$t=0.500$]{
   \includegraphics[scale =0.3] {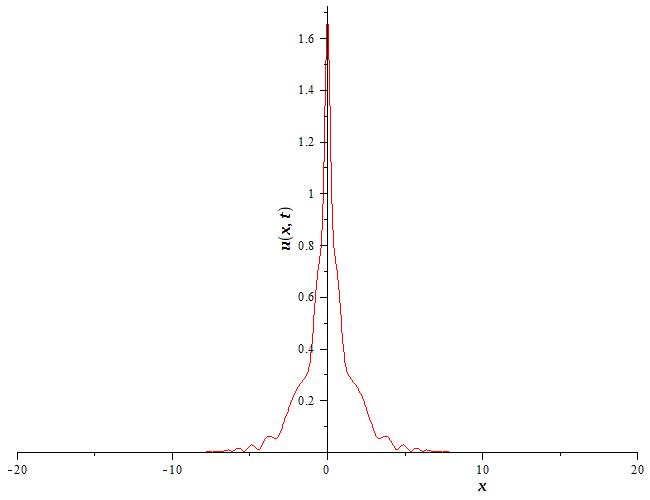}
   \label{bound321}
 }
\caption{Bound state of solitons for $n=3$}
\end{figure}

\noindent
The choice of $n=3$ is studied in the same interval up to the ending time $t=0.6$ using the same discretization parameters. The initial pulse of height $1$ gives birth two new bulges at both sides as the simulation time reaches $t=0.175$ and its height passes over $2$, Fig \ref{bound38}. These two bulges join back to the pulse but then the pulse begins to split in two halves vertically starting in the peak as time proceeds to $t=0.250$, Fig \ref{bound311}. Following the split of the pulse, two bulges begin to appear at the sides of these two halves of the pulse at $t=0.300$, Fig \ref{bound313}, and become evident at $t=0.325$, Fig \ref{bound314}. As time proceeds to $0.400$ and $0.500$, the two halves of the pulse and the bulges join together back and become a unique pulse of height over $1.8$, Fig \ref{bound317} - Fig \ref{bound321}. The similar behaviors keep at the remaining time of the simulation.

\begin{figure}[H]
\centering
   \subfigure[$t=0.175$]{
   \includegraphics[scale =0.3] {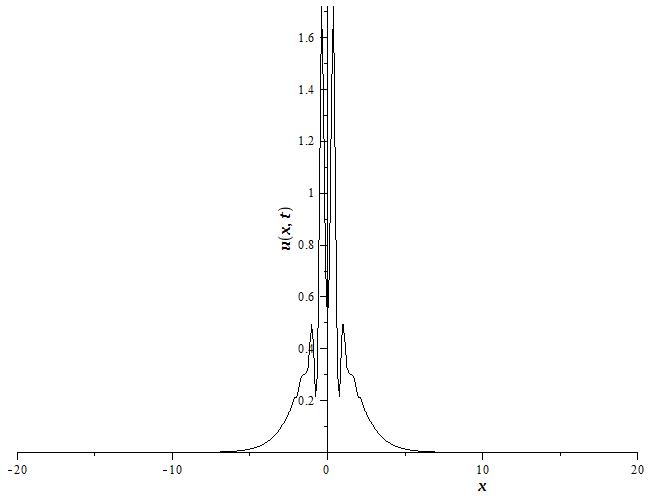}
   \label{bound47}
 }
   \subfigure[$t=0.225$]{
   \includegraphics[scale =0.3] {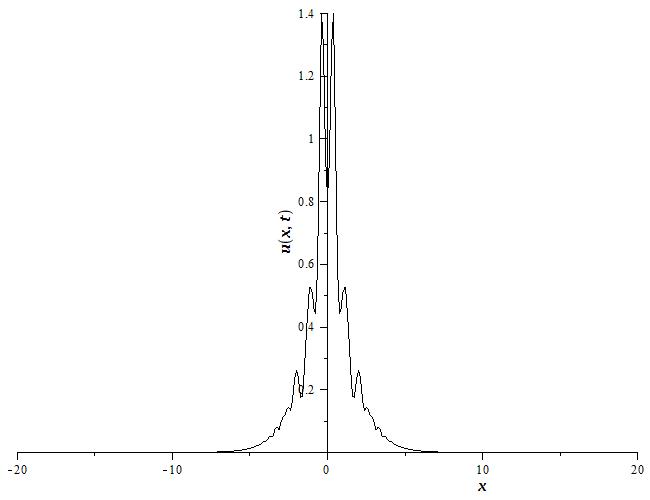}
   \label{bound49}
 }
    \subfigure[$t=0.375$]{
   \includegraphics[scale =0.3] {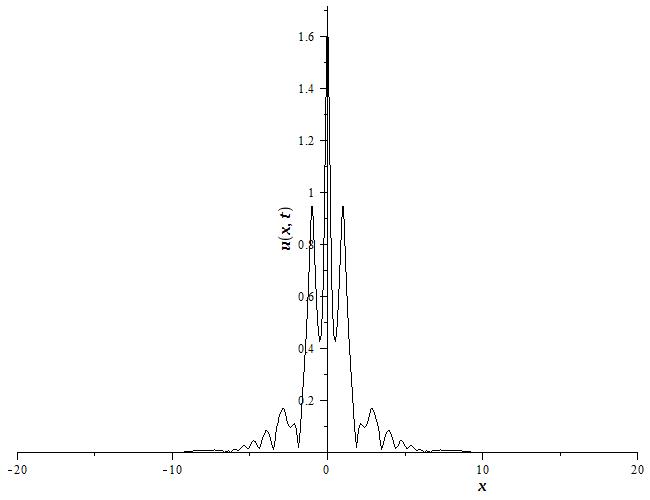}
   \label{bound415}
 }
    \subfigure[$t=0.475$]{
   \includegraphics[scale =0.3] {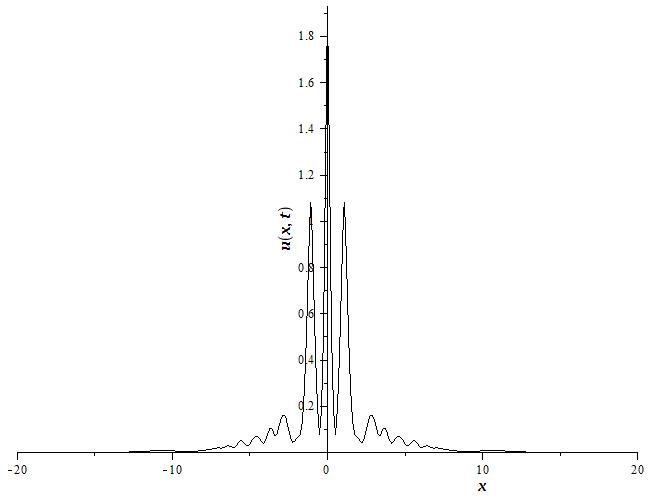}
   \label{bound419}
 }
    \subfigure[$t=0.550$]{
   \includegraphics[scale =0.3] {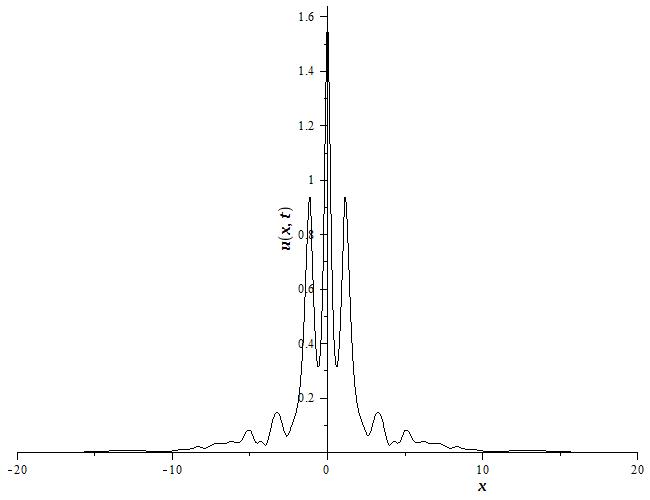}
   \label{bound422}
 }
    \subfigure[$t=0.600$]{
   \includegraphics[scale =0.3] {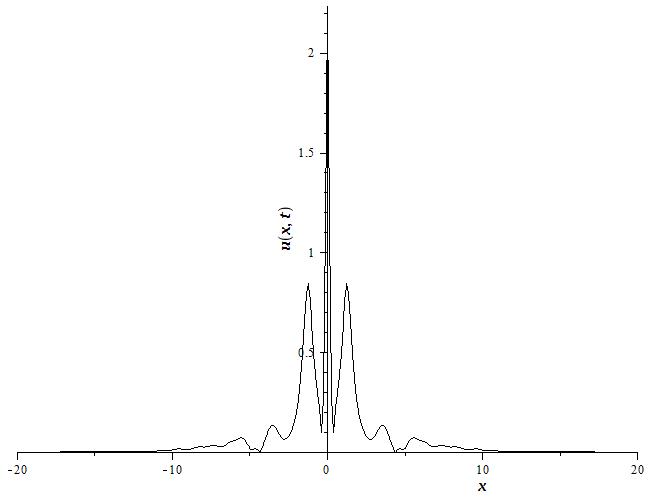}
   \label{bound424}
 }
\caption{Bound state of solitons for $n=4$}
\end{figure}

\noindent
When $n=4$, the initial pulse of unit height splits in two high pulses of heights over $1.6$ with one bulge at the other sides of each at $t=0.175$, Fig \ref{bound47}. As time proceeds the number of bulges increases at both sides of twin longer pulses at time $t=0.225$, Fig \ref{bound49}. At $t=0.375$, two well shaped pulses and a longer pulse of height over $1.6$ between them are observed clearly, Fig \ref{bound415}. There are three smaller bulges at both sides of these three pulses at this time. The twin pulses and the longer pulse positioned between them separate clearly from each other at $t=0.475$, Fig \ref{bound419}. The other smaller pulses propagate along the horizontal axis and get far away from these three longer pulses. The twins and the longer pulse joins together again at $t=0.550$, Fig \ref{bound422}. The heights of the twins approaches $1$ as the longer one gets shorter at this time. When the time is $0.600$, the twin pulses gets smaller in height as the longer one gets higher, Fig \ref{bound424}.

\noindent
In conclusion, when the number $n$ increases, the formation of new solitons are more rapidly and the number of those solitons increases. The solitons do not exceed the imaginary bounds while they are forming and disappearing during the simulations. This status corresponds to the theoretical aspects and the earlier findings reported in the numerical studies.

\noindent
The absolute relative changes of the conservation laws for the experiments obtained for various values of $n$ are illustrated in Table \ref{t5}. The initial values of the conservation laws are determined as $C_1=2$ and $C_3=(2/3)(1-\kappa)$. Both conservation laws remain almost constant for all cases and the absolute relative changes are the indicators of satisfactory as tabulated in the table.

\begin{figure}[H]
 \centering
   \subfigure[$\kappa =2$]{
   \includegraphics[scale =0.3] {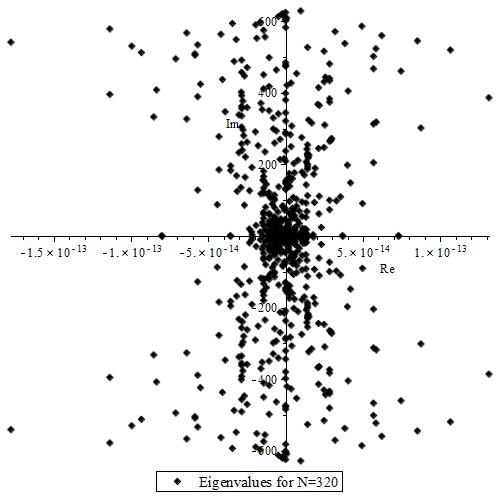}
   \label{EV52}
 }
 \subfigure[$\kappa =3$]{
   \includegraphics[scale =0.3] {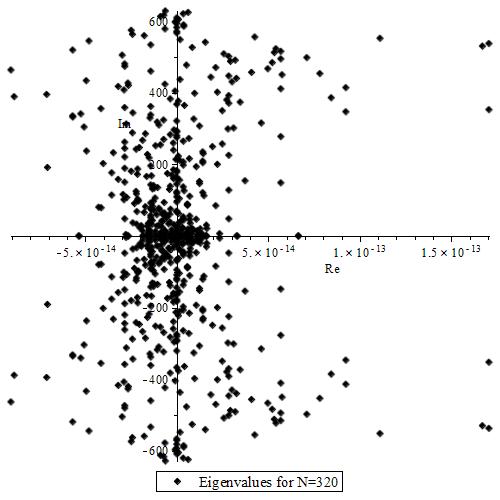}
   \label{EV53}
 }
 \subfigure[$\kappa =4$]{
   \includegraphics[scale =0.3] {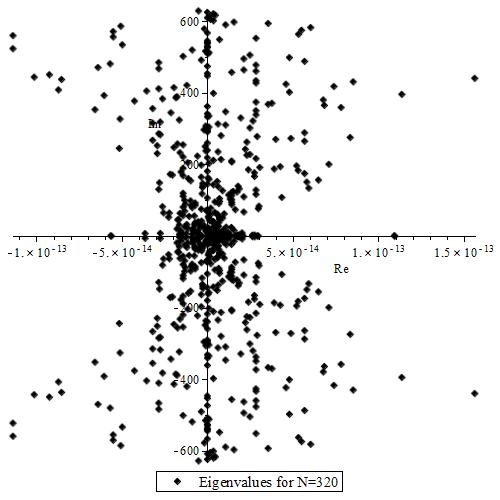}
   \label{EV54}
 }
 \caption{The eigenvalue distribution for bound state of solitons}
\end{figure}

\noindent
The eigenvalue distributions of all cases $\kappa =2,3,4$ are depicted in Fig \ref{EV52}-\ref{EV54}. The maximum and minimum complex components of all eigenvalues are determined as $\pm 628.749$ for all choices of $\kappa$. $\Delta t=0.001$ is sufficient to provide the stability for all methods used to solve the bound state of solitons problem.


\begin{table}[H]
\caption{Absolute relative changes of the conservation laws at $t=0.6$ for the birth of propagating pulse}
\begin{tabular}{lrrrrrrr}
\hline\hline
   Method  & $\Delta x$ &$\Delta t$ & $n$ &  $C_1(0)$  & $C_3(0)$&  $C(C_1(0.6))$ & $C(C_3(6))$ \\
		\hline \hline
		HEUN &$0.125$&$0.001$& $2$ &  $2$  &  $-4.66666666$ &$1.964 \times 10^{-6}$ &  $1.312 \times 10^{-6}$ \\
		RK2  & & & &  &  &$1.964 \times 10^{-6}$ &  $1.311 \times 10^{-6}$ \\
		RK3  & & & &  &  &$6.540 \times 10^{-7}$ &  $4.586 \times 10^{-6}$ \\
		RK4  & &  & & &  &$< 10^{-10}$ &  $6.428 \times 10^{-10}$ \\
		RKF  & & & &  &  &$5.000 \times 10^{-10}$ &  $8.571 \times 10^{-10}$ \\
		CK  & &  & & &  &$5.000 \times 10^{-10}$ &  $8.571 \times 10^{-10}$ \\
		
		HEUN && & $3$ &  $2$  &  $-11.33333333$ &$1.767 \times 10^{-4}$ &  $1.365 \times 10^{-3}$ \\
		RK2  && & &   &  &$1.767 \times 10^{-4}$ &  $1.365 \times 10^{-3}$ \\
		RK3  && & &   &  &$5.629 \times 10^{-5}$ &  $6.353 \times 10^{-4}$ \\
		RK4  && & &   &  &$3.400 \times 10^{-8}$ &  $7.292 \times 10^{-6}$ \\	
		RKF  && & &   &  &$7.500 \times 10^{-9}$ &  $7.238 \times 10^{-6}$ \\
		CK  & & &&   &  &$7.000 \times 10^{-9}$ &  $7.314 \times 10^{-6}$ \\	
		
		HEUN && & $4$ &  $2$  &  $-20.66666666$ &$5.095 \times 10^{-3}$ &  $6.238 \times 10^{-2}$ \\
		RK2  && & &   &  &$5.095 \times 10^{-3}$ &  $6.238 \times 10^{-2}$ \\
		
		RK3  && & &   &  &$7.151 \times 10^{-4}$ &  $8.245 \times 10^{-3}$ \\
		RK4  && & &   &  &$1.587 \times 10^{-6}$ &  $6.432 \times 10^{-4}$ \\	
		
		RKF  & & &&   &  &$5.142 \times 10^{-8}$ &  $6.481 \times 10^{-4}$ \\
		CK  & & &&   &  &$7.675 \times 10^{-7}$ &  $6.370 \times 10^{-4}$ \\		
		\hline
CDQ\cite{korkmaz1} & $0.125$ &  $0.0025$ & $3$ & && $3.282\times 10^{-6}$ & $5.201\times 10^{-4}$ \\   
CDQ\cite{korkmaz1} & $0.125$ &  $0.0025$ & $4$ & && $1.922\times 10^{-4}$ & $1.558\times 10^{-3}$ \\   
PDQ\cite{korkmaz2} & $0.1$ &  $0.01$  & $4$ & && $1.174\times 10^{-7}$ & $5.127\times 10^{-5}$ \\ 		
\hline\hline
\end{tabular}
\label{t5}
\end{table}

\section{Conclusion}
\noindent
Some dynamic problems constructed on the nonlinear cubic NLS equation are solved by differential quadrature method based on sine cardinal functions. The dimension reduced form of the equation is integrated using various methods of different orders. The absolute relative changes are computed in each case to validate the accuracy of the method even analytical solutions do not exist. 

\noindent
The discretization parameter $\Delta t$ is successfully selected to satisfy the stability condition. The instability of some lower order methods, particularly, explained by coinciding the theoretical aspects and matrix stability analysis via eigenvalues for all experiments.

\noindent
The comparison with earlier studies covering some differential quadrature techniques indicate that the proposed algorithms also generate acceptable results. In many cases, the results seem more accurate even when less grids are used. 

\noindent
When compared with the collocation or finite element methods, the programming is easier owing to the main logic that approximates directly to the derivative. Particularly, it also enables using higher order time integration techniques more easily. On the other hand, memory allocation and number of algebraic calculations can be discussed when compared the other method families.

\section{Acknowledgment} This study (Research Project Number BAP-FF160316B19) is supported by The Scientific Research Projects Unit at Çankırı Karatekin University, Turkey.

\end{document}